\documentclass[final]{amsart}
\usepackage{geometry}                % See geometry.pdf to learn the layout options. There are lots.
\usepackage{csquotes}
\usepackage[english]{babel}
\usepackage{amssymb}
\usepackage{amsmath}
\usepackage{amsthm}
\usepackage[pdfstartview=Fit]{hyperref}
\usepackage[color]{showkeys}

\newcommand \FF {\mathcal{F}}
\newcommand \GG {\mathcal{G}}
\newcommand \HH {\mathcal{H}}
\newcommand \EE {\mathcal{E}}
\newcommand \XX {\mathcal{X}}

\DeclareMathOperator{\Aut}{Aut}
\DeclareMathOperator{\Hom}{Hom}
\DeclareMathOperator{\Out}{Out}
\DeclareMathOperator{\Inn}{Inn}
\DeclareMathOperator{\Mor}{Mor}
\DeclareMathOperator{\Syl}{Syl}

\newtheorem{defn}{Definition}[section]
\newtheorem{thm}[defn]{Theorem}

\newtheorem{lemma}[defn]{Lemma}
\newtheorem{prop}[defn]{Proposition}

\newtheorem{claim}{Claim}[defn]

\newtheorem{notation}[defn]{Notation}

\newenvironment{claimproof}[1][Proof of claim]
  {%
    \proof[#1]%
  }
  {%
    \endproof%
  }

\definecolor{refkey}{rgb}{1,0,1}
\definecolor{labelkey}{rgb}{0,0,1}

\begin{document}
\title{Fusion systems on a Sylow \lowercase{$p$}-subgroup of $SU_4(\lowercase{p})$}
\author{Raul Moragues Moncho}
\address{School of Mathematics, University of Birmingham, Edgbaston, Birmingham B15 2TT, United Kingdom}	
%\curraddr{Calle Pintor Cabrera, n\'{u}mero 15, 6 izquierda, 03003, Alicante, Spain}
\email{raul.moragues@gmail.com}
%\urladdr{...} %optional
%\dedicatory{To Chris Parker, for everything} %optional
%\date{\today}
%\thanks{...} %optional
%\translator{...} %Ñ3
%\subjclass[2010]{20D20,  20D05}

\begin{abstract}
We determine, for $p$ odd, all saturated fusion systems $\mathcal{F}$ on a Sylow $p$-subgroup $S$ of the unitary group $SU_4(p)$ and we prove that they are all realizable by finite groups. In particular, we prove that $S$ does not support any exotic fusion systems.
 \end{abstract}
\maketitle

\section{Introduction}
Let $S$ be a Sylow $p$-subgroup of $SU_4(p)$. In this paper we classify the saturated fusion systems $\FF$ on $S$. % satisfying $O_p(\FF)=1$. 
We assume 
%this hypothesis and adopt this notation , as well as 
some basic background about fusion systems throughout the paper. 
Our basic reference for fusion systems is \cite{AKO}. %The case where $p=3$ has been studied by Baccanelli, Franchi and Mainardis in \cite{BacMcL}.  
The case $p=2$ with $O_2(\FF) = 1$ is resolved as part of \cite[Proposition 5.1]{OliverSectionalRank4}, as $UT_4(2)$, a Sylow $2$-subgroup of $SL_4(2)$, is isomorphic to a Sylow $2$-subgroup of $SU_4(2)$ by a straightforward computation. The case $p=3$ with $O_3(\FF) = 1$ has been completed in \cite{BacMcL}. Hence, in most of this paper, we restrict our attention to the case when $p\geq 5$.

The group $S$ contains an extraspecial subgroup of index $p$. In particular, it is case (2) of \cite[Main Theorem]{RMMThesis} and contributes to the project of classifying saturated fusion systems on $p$-groups containing an extraspecial subgroup of index $p$, which was the object of study of my PhD Thesis. The result of this paper will be applied in forthcoming work on this topic.

It will also form part of the classification of fusion systems of sectional $p$-rank $4$ for odd primes $p$. This case is interesting because there is an $\FF$-essential subgroup $V$ such that $|N_S(V)/V| = p^2$, which does not happen in the cases with smaller sectional $p$-rank.
Since there are two $\FF$-essential subgroups, it can be seen as a contribution towards the classification of rank 2 fusion systems on Sylow $p$-subgroups of groups of Lie type, studied in %such as those of $G_2(p)$ studied in 
\cite{ParkerSemeraroG2} and  %$GL_3(p)$ classified in 
\cite{RV}. 

All of these contributions add to our knowledge of saturated fusion systems defined on $p$-groups for odd primes $p$ and so extend our understanding of how exotic fusion systems arise at odd primes [AKO, Problem III.7.4] and [AO, Problem 7.6]. 
In the case under consideration, we prove that no exotic fusion systems arise.

The main theorem of this paper is the following:

\begin{thm}\label{mainthm}
Assume $S$ is a Sylow $p$-subgroup of $SU_4(p)$. Let $\FF$ be a saturated fusion system on $S$ and, if $p=2$, assume further that $O_2(\FF) = 1$. 
%
%OLD Then $\FF$ is known and realizable by a finite group.
%
%More precise: 
Then $\FF$ is realizable by a finite group and either $\FF$ is constrained or $O_p(\FF)=1$. If $O_p(\FF) = 1$ then one of the following holds:
\begin{itemize}
\item$p=2$ and $\FF$ is isomorphic to the fusion system of $GL_4(2)$ or to that of $PSp_4(q)$ for each $q\equiv \pm 3\pmod 8$;

\item $p=3$ and 
%(BFM Table 2 ? ) %. no exotics. 
%	or using their statement: 
%	either constrained ($O_3(\FF)\neq 1$) or 
	$\FF$ is isomorphic to a fusion system $\FF_S(G)$, where $G$ is one of the following:
	\begin{enumerate}
	\item[(i)]  $\tilde G\leq G\leq \Aut(\tilde G)$ with $\tilde G\in \{Mc, U_4(3), Co_2\}$;
	\item[(ii)]  $G = L_6(q)$ with $q\equiv 4,7\pmod 9$, or $G=U_6(q)$ with $q\equiv 2,5\pmod 9$;
	\item[(iii)]  $G = L_6(q)\langle\phi\rangle$ with $q\equiv 4,7\pmod 9$, or $G=U_6(q)\langle\phi\rangle$ with $q\equiv 2,5\pmod 9$, where $\phi$ is a field automorphism of order $2$;
	\end{enumerate}
		Moreover, all the fusion systems in (ii) (respectively in (iii)) realize isomorphic fusion systems.

\item $p\geq5$ and there is a one-to-one correspondence between saturated fusion systems $\FF$ on $S$ 
%with $O_p(\FF)=1$ 
and groups $G$ with $SU_4(p)\leq G \leq \Aut(SU_4(p))$ which realize them. 
\end{itemize}
\end{thm}

In our argument, we consider separately the cases according to whether $\FF$ contains any non-trivial normal $p$-subgroup, that is, whether $O_p(\FF) = 1$ or not. 
When $p$ is odd and $O_p(\FF)\neq 1$,  we will show in Proposition \ref{OpNotTrivial} that there is at most one $\FF$-essential subgroup, and it is normal in $\FF$, whence the Model Theorem implies that $\FF$ is realizable.  The case when $p=2$ and  $O_2(\FF)\neq 1$ has been checked computationally using the algorithms in \cite{ParkerSemeraroAlgorithms}.

When $O_p(\FF) = 1$, the situation is more complicated.
When $p=2$, we note that $PSU_4(2)\cong PSp_4(3)$, and \cite[Proposition 5.1]{OliverSectionalRank4} proves that, if $O_2(\FF) = 1$, then $\FF$ can arise only as the fusion system of $GL_4(2)$, or that of $PSp_4(q)$ for $q\equiv \pm 3\pmod{8}$, and there are only 2 saturated fusion systems $\FF$ on $S$ with $O_2(\FF) =1$.

When $p=3$, \cite[Theorem 1]{BacMcL} implies that no exotic fusion systems arise. This case requires different arguments as our classification for $p\geq 5$, 
and there is more than one possibility for 
$\FF_0$, the smallest fusion system on $S$ which contains $O^{p'}(\Aut_\FF(P))$ for each $P\leq S$, since 
the shape of $O^{p'}(\Aut_\FF(P))$ for some subgoups $P\leq S$ is not uniquely determined.  This gives rise to more realizable saturated fusion systems $\FF$, such as those arising from $L_6(q)$, $U_6(q)$ for appropriate $q$ coprime to $p$, and the sporadic finite simple groups $McL$ and $Co_2$.  There are 
$13$ such fusion systems when $p=3$. We will prove that, when $p\geq 5$, the number of saturated fusion systems $\FF$ with $O_p(\FF) = 1$ is either $5$ when $p\equiv 1\pmod 4$ or $8$ when $p\equiv 3\pmod 4$.

We remark that, as written, the proof depends on the Classification of Finite Simple Groups, which is used in Proposition \ref{speinGL4} to determine the shape of a subgroup of $GL_4(p)$. 

The paper is structured as follows: we begin by introducing the background definitions and results needed in Section \ref{sec:background}.
In Section \ref{sec:structure}, we describe the $p$-group $S$ and some of its subgroups, such as the unique extraspecial subgroup $Q$ of index $p$ and the unique abelian subgroup $V$ of order $p^4$ in $S$. We also consider the automorphism groups of $S$, $Q$ and $V$, and prove results about the action of the relevant subgroups on $Q/Z(Q)$ and $V$.

In Section \ref{sec:essential} we determine that the only subgroups of $S$ which can be $\FF$-essential are $V$ and $Q$, use the results from the previous section to determine $O^{p'}(\Aut_\FF(Q))$ and $O^{p'}(\Aut_\FF(V))$ when they are $\FF$-essential, and then prove that both $V$ and $Q$ are $\FF$-essential exactly when $O_p(\FF) = 1$. 

We begin Section \ref{sec:classification} by concluding the case where $O_p(\FF) \neq 1$ in Proposition \ref{OpNotTrivial}. In the case $O_p(\FF) =1$, the following is our main result.

\begin{thm}\label{CaseU} Assume $p\geq 5$ and $S$ is a Sylow $p$-subgroup of $SU_4(p)$. Then there is a one-to-one correspondence between saturated fusion systems $\FF$ on $S$ with $O_p(\FF)=1$ and groups $G$ with $SU_4(p)\leq G \leq \Aut(SU_4(p))$ which realize them. In particular, there are no exotic fusion systems $\FF$ on $S$ with $O_p(\FF)=1$.
\end{thm}

We begin the proof of Theorem \ref{CaseU} by considering the subgroups of $\Aut_\FF(S)$ generated by those elements which restrict to maps in one of $O^{p'}(\Aut_\FF(V))$ or $O^{p'}(\Aut_\FF(Q))$, focusing on their $p'$-parts, and determine their actions on various sections of $S$. We then consider $\FF_0$, the smallest saturated fusion system on $S$ with $O_p(\FF_0) = 1$, which exists by the Alperin-Goldschmidt Fusion Theorem. We use the assumption that $p\geq 5$ and the subgroups considered above to determine $\Aut_{\FF_0}(S)$ uniquely as a subgroup of $\Aut_{\Aut(SU_4(p))}(S)$. With a fixed copy of $\Aut_{\FF_0}(S)$, we consider the subgroup of morphisms centralizing both $Z$ and $S/Q$ to determine $O^{p'}(\Aut_{\FF_0}(V))$ and $O^{p'}(\Aut_{\FF_0}(Q))$ uniquely, thus establishing the uniqueness of $\FF_0$ up to isomorphism. Finally, we realize $\FF_0$ as the fusion system of $PSU_4(p)$, and the fusion system of $\Aut(PSU_4(p))$ as the largest possible saturated fusion system on $S$, which allows us to realize every saturated fusion system by an intermediate group.

Theorem \ref{CaseU} is proved as Theorem \ref{ClassificationCaseU}. Theorem \ref{mainthm} then follows putting together  \cite[Proposition 5.1]{OliverSectionalRank4}, \cite[Theorem 1]{BacMcL}, Proposition \ref{OpNotTrivial} and Theorem \ref{CaseU}.

We remark that we write maps on the right side of the argument.

\section{Preliminaries}\label{sec:background}

Given a finite group $G$, the fusion system $\FF_S(G)$ of $G$ on $S$ is a category which encodes the $p$-local structure of $G$ on one of its Sylow $p$-subgroups $S$.  The objects of $\FF_S(G)$ are the subgroups of $S$ and, given $P,Q\leq S$, the morphisms of $\FF_S(G)$ are given by the conjugation maps $c_g$ induced by $g\in G$ with domain $P$ and codomain $Q$, that is

$$\Mor_{\FF_S(G)}(P,Q) = \{c_g\mid g\in G \text{ with } P^g\leq Q\}.$$
$\FF_S(G)$ is a foundational example of a \emph{fusion system} on $S$, as defined in \cite[Definition 2.1]{AKO}. 

A proper subgroup $H < G$ of a finite group $G$ is \emph{strongly $p$-embedded} in $G$ if $p$ divides $|H|$ and, for each $g\in G\setminus H$, $p$ does not divide $|H\cap H^g|$.

We now introduce the basic concepts of the theory of fusion systems which we will use.

\begin{defn}\label{bigdef}Let $\FF$ be a fusion system on a $p$-group $S$, and $P,Q\leq S$. Then,
\begin{enumerate}
\item $P$ is \emph{fully $\FF$-automized} if $\Aut_S(P)\in \Syl_p(\Aut_\FF(P))$;
\item $P^{\FF} := \{P\alpha \mid \alpha\in \Hom_\FF(P,S)\}$;
\item $P$ is \emph{fully $\FF$-normalized} if $|N_S(P)|\geq |N_S(Q)|$ for all $Q\in P^\FF$;
\item $P$ is \emph{fully $\FF$-centralized} if $|C_S(P)|\geq |C_S(Q)|$ for all $Q\in P^\FF$;
\item $P$ is \emph{$\FF$-centric} if $C_S(Q) = Z(Q)$ for all $Q\in P^\FF$;
\item $P$ is \emph{$\FF$-essential} if $P$ is $\FF$-centric and fully $\FF$-normalized, and $\Out_\FF(P)$ contains a strongly $p$-embedded subgroup;
\item if $\alpha\in \Hom_\FF(P,Q)$ is an isomorphism, then % \emph{$\alpha$-extension control subgroup} of $S$ is 
$N_\alpha = \{g\in N_S(P)\mid \alpha^{-1}c_g\alpha\in \Aut_S(Q)\};$

\item $Q$ is \emph{$\FF$-receptive} if for all isomorphisms $\alpha\in \Hom_\FF(P,Q)$ there exists $\widetilde\alpha\in \Hom_\FF(N_\alpha, S)$ such that $\widetilde \alpha|_P=\alpha$.
\end{enumerate}
\end{defn}

We can now define a \emph{saturated fusion system}, of which $\FF_S(G)$ defined earlier is an example. 

\begin{defn}

A fusion system $\FF$ is \emph{saturated} if for each subgroup $P\leq S$, there exists $Q\in P^\FF$ such that
$Q$ is both fully $\FF$-automized and 
 $\FF$-receptive.
\end{defn}

Those saturated fusion systems that arise as $\FF_S(G)$ for a finite group $G$ with Sylow $p$-subgroup $S$ are called \emph{realizable}, and those that do not are called \emph{exotic}.
From now on, we restrict our attention to saturated fusion systems. We now present some of their properties. 
By \cite[Lemma 2.6~(c)]{AKO}, a subgroup $Q$ of $S$ is fully $\FF$-normalized if and only if it is fully $\FF$-automized and $\FF$-receptive. In particular, $\FF$-essential subgroups are both fully $\FF$-automized and $\FF$-receptive. The following lemma is a consequence of $\FF$ being saturated.

\begin{lemma}\label{essentiallift}

Let $\FF$ be a saturated fusion system on the $p$-group $S$. Suppose that $E\leq S$ is fully $\FF$-normalized. Then  every element of $N_{\Aut_\FF(E)}(\Aut_S(E))$ lifts to an element of $\Aut_\FF(N_S(E))$.
\end{lemma}	
\begin{proof} Follows directly from the property that a fully
$\FF$-normalized subgroup is receptive.
%This is a consequence of the of the surjectivity property, see for example \cite[Definition 6.3]{Craven} and the discussion following it.
\end{proof}

The starting point in the classification of saturated fusion systems is the Alperin-Goldschmidt Fusion Theorem. In order to state it, we first define what it means for certain maps to generate a fusion system.
 If $X$ is a set of injective morphisms between various subgroups of $S$, then we define $\langle X \rangle$ to be the fusion system obtained by intersecting all the fusion systems on $S$ which have the members of $X$ as morphisms.

\begin{thm}[{Alperin-Goldschmidt Fusion Theorem \cite[Theorem I.3.5]{AKO}}] \label{Alperin}
 Suppose $\FF$ is a saturated fusion system on a $p$-group $S$. Then 
$$\FF = \langle \Aut_\FF(S), \Aut_\FF(E) \mid E \text{ is $\FF$-essential}\rangle .$$
\end{thm}

The Alperin-Goldschmidt Fusion Theorem allows us to focus on the $\FF$-essential subgroups. In particular, we will repeatedly use the fact that  if $E$ is $\FF$-essential then, as $\Out_\FF(E)$ has a strongly $p$-embedded subgroup, we have $O_p(\Out_\FF(E)) = 1$, which implies that $O_p(\Aut_\FF(E)) = \Inn(E)$. The Alperin-Goldschmidt Fusion Theorem Theorem can be refined using Frattini's Argument and Lemma \ref{essentiallift} to obtain the following.

\begin{lemma}[{\cite[Lemma 3.4]{BCGLO2}}]\label{Op'AndSGenerate}\leavevmode If $\FF$ is saturated then $$\FF = \langle  \Aut_\FF(S),  O^{p'}(\Aut_\FF(E))\mid E\text{ is } \FF\text{-essential}\rangle.$$
\end{lemma}

We now turn to normality in fusion systems.

\begin{defn} Let $\FF$ be a fusion system on $S$, and $Q\leq S$. Then:

\begin{itemize}
\item $Q$ is \emph{normal} in $\FF$, denoted by $Q\trianglelefteq \FF$, if $Q\trianglelefteq S$ and, for all $P,R\leq S$ and all $\phi\in \Hom_\FF(P,R)$, $\phi$ extends to a morphism $\overline\phi\in \Hom_\FF(PQ,RQ)$ such that $Q\overline\phi = Q$;

\item $O_p(\FF)\trianglelefteq S$ denotes the largest subgroup of $S$ which is normal in $\FF$.
\end{itemize}
\end{defn}

We will determine that $O_p(\FF) = 1$ by using the following result.

\begin{prop}[{\cite[Proposition I.4.5]{AKO}}]\label{AKONormal} Let $\FF$ be a saturated fusion system on $S$. Then for any $Q\leq S$, $Q\trianglelefteq \FF$ 
if and only if 
for each $P\leq S$ which is $\FF$-essential or equal to $S$, we have that $Q\leq P$ and $Q$ is $\Aut_\FF(P)$-invariant.
\end{prop}

\begin{defn} \label{AutSFromE} Suppose $\FF$ is saturated. We define the following.
\begin{enumerate}
\item If $E$ is $\FF$-essential  
then 
$\Aut_\FF^E(S) := \langle \alpha \in \Aut_\FF(S) \mid \alpha|_E\in O^{p'}(\Aut_\FF(E))\rangle$;

\item $\Aut_\FF^0(S) := \langle \Aut_\FF^E(S), ~\Inn(S) \mid E\text{ is } \FF\text{-essential}\rangle$;

\item $\FF_0 := \langle O^{p'}(\Aut_\FF(E)), \Aut_\FF^0(S) \mid E\text{ is } \FF\text{-essential} \rangle\subseteq \FF$; 

\item $\Gamma_{p'}(\FF) := \Aut_\FF(S)/\Aut_\FF^0(S)$;

\item Let $\EE$ be a subcategory of $\FF$ which is itself a saturated fusion system over a subgroup $T\leq S$.
Then $\EE$ has \emph{index prime to $p$} in $\FF$ if $T=S$ and $\Aut_\EE(P)\geq O^{p'}(\Aut_\FF(P))$ for each $P\leq S$.
\end{enumerate}
\end{defn}
Thus $\Aut_\FF^E(S)\leq \Aut_\FF^0(S)$ is the subgroup of automorphisms that are contributed to $\Aut_\FF(S)$ by $O^{p'}(\Aut_\FF(E))$.  We remark that $\Aut_{\FF_0}(S) = \Aut_\FF^0(S)$.
The following allows us to characterize subsystems of index prime to $p$.
\begin{thm}[{\cite[Theorem I.7.7]{AKO}}]\label{PrimeIndexFusion}Suppose that $\FF$ is a saturated fusion system on $S$. Then there is a one-to-one correspondence between saturated fusion subsystems of $\FF$ on $S$ of index prime to $p$ and subgroups of $\Gamma_{p'}(\FF)$. In particular, there is a unique minimal saturated fusion subsystem of index prime to $p$, which we denote by $O^{p'}(\FF)$. We have $\Aut_{O^{p'}(\FF)}(S) = \Aut_\FF^0(S)$.
\end{thm}

Note that 
by definition and the Alperin-Goldschmidt Fusion Theorem,
 $\FF_0$ is the smallest fusion system on $S$ which contains $O^{p'}(\Aut_\FF(P))$ for each $P\leq S$, 
 and $\FF_0 = O^{p'}(\FF)$. 
We have $\FF = \langle \FF_0, \Aut_\FF(S)\rangle$, so we will construct the fusion system $\FF_0$, show it is saturated by finding a group realizing it, determine the largest possible candidate for $\Aut_\FF(S)$, and then use Theorem \ref{PrimeIndexFusion}  to obtain all subsystems of $p'$-index as intermediate fusion systems.

The relationship between $\Aut_\FF^E(S)$ and $O^{p'}(\Aut_\FF(E))$ is made clear by the following result. %, which uses an equivalent formulation of saturation via the surjectivity property (see \cite[\textsection 6.1]{Craven} for details).

\begin{lemma} %[{\cite[Lemma 2.37]{RMMThesis}}]
\label{OutFSBijection}
If $\FF$ is a saturated fusion system on $S$ and  $E\trianglelefteq S$ is $\FF$-centric and normalized by $\Aut_\FF(S)$, then there are isomorphisms
\begin{align*} &\Aut_\FF(S)/C_{\Inn(S)}(E)\cong N_{\Aut_\FF(E)}(\Aut_S(E)), \\
&\Out_\FF(S)\cong N_{\Aut_\FF(E)}(\Aut_S(E))/\Aut_S(E) \text{, and }\\
 &\Out_\FF^E(S)\cong N_{O^{p'}(\Aut_\FF(E))}(\Aut_S(E))/\Aut_S(E).\end{align*}
\end{lemma}

\begin{proof} It is a consequence of Thompson's $A\times B$-Lemma that $C_{\Aut_\FF(S)}(E)$ does not contain any nontrivial $p$'-subgroup. Hence, $C_{\Aut_\FF(S)}(E)\leq \Inn(S)$ and the assertion follows from Lemma \ref{essentiallift}.
\end{proof}

We now introduce the concepts of morphisms and isomorphisms between fusion systems.

Two fusion systems $\FF$ and $\EE$ on $S$ are isomorphic if there exists $\alpha\in \Aut(S)$ such that for all $P, Q \leq S$,
$$\Hom_\EE(P\alpha,Q\alpha) = \{\alpha^{-1}|_{P\alpha}\theta\alpha \mid \theta\in \Hom_\FF(P,Q)\}.$$

We note that an isomorphism of fusion systems preserves saturation by \cite[Lemma II.5.4]{AKO}.
An easy isomorphism of fusion systems is given in the following situation.

\begin{lemma}\label{IsomorphismQuotient}
Let $G$ be a finite group, $S\in \Syl_p(G)$. Let $N\trianglelefteq G$ with $p\nmid |N|$, $\bar G = G/N$, $\bar S \in \Syl_p(\bar G)$. Then $\FF_S(G)\cong \FF_{\bar S}(\bar G)$.
\end{lemma}

We now turn our attention to the group theoretic methods that we will use in the paper. We begin with some basic results about coprime action.

\begin{lemma}[{\cite[Corollary 5.3.3]{Gorenstein}}]
\leavevmode\label{Oliver1.7}Fix a prime $p$, a finite $p$-group $S$, and a group $G\leq \Aut(S)$ of automorphisms of $S$. Let $S_0\trianglelefteq S_1\trianglelefteq \dots \trianglelefteq S_m = S$ be a sequence of subgroups, all normal in $S$ and normalized by $G$, such that $S_0\leq \Phi(S)$. Let $H\leq G$ be the subgroup of those $g\in G$ which act via the identity on $S_i/S_{i-1}$ for each $1\leq i\leq m$. 
Then $H$ is a normal $p$-subgroup of $G$. In particular, $C_{\Aut(S)}(S/\Phi(S))$ is a normal $p$-subgroup of $\Aut(S)$.
\end{lemma}

We will  often use this result in order to prove that $O_p(\Out_\FF(E))\neq 1$, which implies that $\Out_\FF(E)$ cannot have a strongly $p$-embedded subgroup, thus proving that various subgroups $E$ of $S$ cannot be $\FF$-essential.

We finish this section with 
a result about groups with a strongly $p$-embedded subgroup whose Sylow $p$-subgroups are not cyclic. 

\begin{thm}\label{almostsimplespe}
Assume that $G$ is a finite group, $H<G$ is strongly $p$-embedded, and $H$ contains an elementary abelian subgroup of order $p^2$. Then $O_{p'}(G)\leq H$,  $H/O_{p'}(G)$ is strongly $p$-embedded in $G/O_{p'}(G)$, and $O^{p'}(G/O_{p'}(G))$ is a nonabelian almost simple group.
\end{thm}

\begin{proof}Assume that $H$ is strongly $p$-embedded in $G$. Let $\overline G = G/O_{p'}(G)$. 
Assume $T$ is a normal $p'$-subgroup of $G$. Let $A\leq H$ be elementary abelian of order $p^2$. Then by \cite[Theorem 6.2.4]{Gorenstein} we have $T = \langle C_T(a)\mid a\in A\setminus 1\rangle$.

As $H$ is strongly $p$-embedded in $G$, we have $T\leq H$ by \cite[Proposition 17.11]{GLS2}.  In particular, $O_{p'}(G)\leq H$. Let $\overline H  = H/O_{p'}(G)$ and  $\overline S = SO_{p'}(G)/O_{p'}(G)$. 
%As $\FF_S(G)\cong \FF_{\overline S}(\overline G)$, it follows that 
%$\overline H$ is strongly $p$-embedded in $\overline G$.
Since $p\mid H$, we have $p\mid \overline H$. Now let  $\overline x = xO_{p'}(G)\in \overline G\setminus \overline H$, then $x\in G\setminus H$ and, since $H$ is strongly $p$-embedded in $G$, we have $p\nmid |H\cap H^x|$. Therefore, $p\nmid |\overline H\cap \overline H^{\overline x}|$, so $\overline H$ is strongly $p$-embedded in $\overline G$.

So now assume $O_{p'}(G) = 1$, fix $S\in \Syl_p(G)$, and let $T$ be a minimal normal subgroup of $G$. As $G$ has a strongly $p$-embedded subgroup $H$, $O_p(G) = 1$ by \cite[Proposition 17.11]{GLS2}. Then $T$ is not a $p$-group and not a $p'$-group either as $O_{p'}(G) = 1$. As $T$ is minimal normal, it is characteristically simple. Thus, as $T$ is nonabelian, $T$ is a direct product of isomorphic nonabelian simple groups $T = L_1\times \ldots \times L_k$ by \cite[Lemma 8.2]{aschbacher}.

If $G$ is not almost simple then either $k>1$ or there is another minimal normal subgroup distinct from $T$. If such a minimal normal subgroup exists, denote it by $T_2$, if there is no such subgroup, let $T_2 = 1$. 
% If $k= 1$ then $T$ is simple.
%If there is another minimal normal subgroup disctinct from $T$, denote it by $T_2$, if there is no such subgroup, let $T_2  = 1$. 
Then, 
%If $k>1$ or there is another minimal normal subgroup $T_2$, 
we have $[T, T_2]\leq T\cap T_2 = 1$, so $T\times T_2 \leq soc(G)\leq G$. Let $S_1 = S\cap L_1\in \Syl_p(L_1)$ and $S_2 = S\cap (L_2\times ... \times L_k\times T_2)$ which is in $\Syl_p(L_2\times ... \times L_k\times T_2)$. As $O_{p'}(G) = 1$, $S_1$ and $S_2$ are nontrivial. Hence $L_1\leq N_G(S_2)$ and $L_2\times ...\times L_k\times T_2\leq N_G(S_1)$. Let $S_0 = S\cap T\in \Syl_p(T)$. Since $H$ is strongly $p$-embedded in $G$ then $N_G(S_i)\leq H$, so  $TT_2= L_1\times L_2\times ... \times L_k\times T_2 \leq H$.  Then  we have, by the Frattini Argument, $N_G(S_0)T = G$, but $N_G(S_0)T\leq H\neq G$, a contradiction. 
Thus $T$ is the unique minimal normal subgroup of $G$, and is a nonabelian simple group. Then, $G$ is almost simple. 
\end{proof}

We conclude this section with two straightforward results about the eigenvalues and eigenspaces of modules. 

\begin{lemma}\label{AuxiliaryMeier1} Let $G$ be a group, $K$ be a field, $S\in \Syl_p(G)$, and $V_1$, $V_2$ be $KG$-modules. 
Let $\theta: V_1\rightarrow V_2$ be a $KG$-module isomorphism. Then the eigenvalues of $r\in N_G(S)$ on $C_{V_1}(S)$ and on $C_{V_2}(S)$ are the same.
\end{lemma}

\begin{proof}
Let $v\in C_{V_1}(S)$, then $(v\theta)s = (vs) \theta = v\theta$ for all $s\in S$, so $C_{V_1}(S)\theta = C_{V_2}(S)$.
Further, let $v$ be a $\lambda$-eigenvector for $r$, then $(v\theta)r = (vr)\theta = (\lambda v)\theta = \lambda (v\theta)$, so $v\theta$ is a $\lambda$-eigenvector for $r$. The same argument with $\theta^{-1}$ gives this property for $\theta^{-1}$, hence the lemma is proved.
\end{proof}

\begin{lemma}\label{AuxiliaryMeier2}Let $H\cong SL_2(p)$ and $S\in \Syl_p(H)$, let $V$ be the natural module for $H$ and let  $r\in N_H(S)\setminus S$. %Then $r$ normalises $S$ and $S^h$ for some $h\in H$ and, if $o(r)\neq 2$, then $C_V(S)$ and $C_V(S^h)$ are eigenspaces for distinct eigenvalues of the action of $r$ on $V$.
Then there exists $h\in H\backslash N_H(S)$ such that $r$ normalizes $S^h$. Moreover, if $o(r)\neq 2$, the subspaces $C_V(S)$ and $C_V(S^h)$ are eigenspaces of $r$ for distinct eigenvalues $\lambda$ and $\lambda^{-1}$.
\end{lemma}
\begin{proof} By Sylow's Theorems, we may assume that $S$ consists of unipotent lower triangular matrices. Then, as $N_H(S)\cong C_p\rtimes C_{p-1}$ contains a unique conjugacy class of complements to $S$ and $r\in N_H(S)\setminus S$, we may assume $r = \left(\begin{smallmatrix}\lambda&0\\0&\lambda^{-1} \end{smallmatrix}\right)$ for some $\lambda\in GF(p)$, hence $r$ normalises the subgroup of unipotent upper triangular matrices, which is another Sylow $p$-subgroup of $H$ and hence is $S^h$ for some $h\in H\setminus N_H(S)$. Further, $r$ has eigenspaces $C_V(S)$ and $C_V(S^h)$ with respective eigenvalues  $\lambda$ and $\lambda^{-1}$, which are distinct unless $\lambda^2=1$, that is $o(r) =2$. 
%Since $H\cong SL_2(p)$ acts  $2$-transitively on its Sylow $p$-subgroups, we may assume that $S$ and $S^h$ consist of unipotent lower and upper triangular matrices respectively. 
\end{proof}

\section{Structure of \texorpdfstring{$S$}{\textit{S}}}\label{sec:structure}

Throughout the rest of this paper we assume that $p\geq 5$ unless explicitly stated otherwise.
The unitary group $SU_4(p)$ has Sylow $p$-subgroups $S$ of order $p^6$ by \cite[Theorem 2.2.9]{GLS3} and, by \cite[Table 8.10]{BHRD}, contains maximal parabolic subgroups containing $S$ of shapes 

$$P_1\cong p^{1+4}_+:SU_2(p) : (p^2-1) \text{ and } P_2\cong C_p^4:SL_2(p^2):(p - 1).$$ 

We denote by $Q := O_p(P_1)\cong p^{1+4}_+$ and $V := O_p(P_2)\cong C_p^4$.

Further, by \cite[Proposition 3.3.1]{GLS3}, $S$ has nilpotency class $3$, and its upper and lower central series coincide, with $S' = Z_2(S)$ of order $p^3$ and $Z := Z(S) = [S,S,S]$ of order $p$. As we assume $p \geq 5$, and $SU_4(p)$ embeds into $GL_4(p^2)\leq GL_p(p^2)$, which has Sylow $p$-subgroups of exponent $p$, $S$ has exponent $p$.

From these details it can be seen that $S$ is a semidirect product of $Q$ by a group of order $p$ whose generator acts on the symplectic vector space $Q/Z$ with Jordan form of two blocks of size $2$. We can also describe $S$ as $V\rtimes T$ where $T$ is isomorphic to a Sylow $p$-subgroup of $\Omega_4^-(p)\cong PSL_2(p^2)$. This explains parts (1)-(3) of the following lemma.

\begin{lemma}\label{LUPropertiesOfS} The following hold:
\begin{enumerate}
\item The order of $S$ is $p^6$.
\item $S$ has nilpotency class $3$ and the terms of the upper and lower central series of $S$ are $Z := Z(S) = [S,S,S]$ of order $p$ and 
$Z_2(S)=\Phi(S)=S'$ of order $p^3$. 
\item $S$ has exponent $p$.
\item $Q$ is the unique maximal subgroup of $S$ with $|Q'| = p$. In particular, $Q$ is the unique  extraspecial subgroup of index $p$ in $S$, and $Q$ is characteristic in $S$.
\item $V$ is the unique abelian subgroup of order $p^4$ in $S$. In particular, $V$ is characteristic in $S$.
\item Let $\XX := \{ M\leq S\mid V<M<S\}$. Then $|\XX| = p+1$,  every $M\in \XX$ satisfies $|Z(M)| = |M'| = p^2$. 
Further, there is an element of $\Aut(S)$ of order $p+1$ which permutes transitively the elements of $\XX$. In particular, the elements of $\XX$ are all isomorphic and have $Z(M)\neq M'$. 
\item Let $M\leq S$ be a maximal subgroup with $M\notin \XX$. Then $Z(M) = Z(S)$.
\end{enumerate}

\end{lemma}

\begin{proof}Parts (1), (2) and (3) are proved above.  Part (4) follows from part (2), as if there was another maximal subgroup $M$ of $S$ with $|M'|=p$, then $M'=Z$, and $Z(S/Z)\geq Q/Z\cap M/Z$, which has order $p^3$. This contradicts $|Z_2(S)| = p^3$.

Now if $W\neq V$ was another abelian subgroup of $S$ of order $p^4$ then $V\cap W\leq Z(VW)$. Hence, as $|Z(S)|=p$, we would have $VW<S$, which would imply $Z(S/Z)\geq Q/Z\cap VW/Z$, again contradicting $|Z_2(S)| = p^3$, hence part (5) holds.

For part (6), we consider $P_2\cong C_p^4 : SL_2(p^2) : C_{p-1}$, and the corresponding maximal subgroup $M$ of $\Aut(SU_4(p))$, where there is an element $\theta$ of $SL_2(p^2)$ of order $p^2-1$ normalizing $S$. Then, $\theta$ normalizes $V$ and acts on a complement to $V$ in $S$. Furthermore, $\theta^{p-1}$ acts transitively on the $p+1$ maximal subgroups of $S$ containing $V$, whence all elements of $\XX$ are isomorphic. Note that the central involution of the $SL_2(p^2)$ acting on $V$ is the central involution of $SU_4(p)$, and as $PSL_2(p^2)\cong \Omega_4^-(p)$, we see that the Jordan form of the action on $V$ of every $p$-element in $P_2$ not in $V$ has one Jordan block of size $3$ and one of size $1$  
by \cite[Theorem 3.1~(ii)]{LiebeckSeitzUnipotent}.
Hence, $Z(M)$ and $M'$ do not coincide and both have order $p^2$.

Finally, we turn to part (7). 
Since $S'\leq M \leq C_S(Z(M))$, we see that $Z(M)\leq C_S(S') = V$, hence $C_S(Z(M))\geq VM = S$, so $Z(M) = Z(S)$. 
\end{proof}

By the Alperin-Goldschmidt Fusion Theorem \ref{Alperin}, in order to determine the fusion systems $\FF$, we will have to determine the automorphism groups of the $\FF$-essential subgroups and $S$. Hence we now describe the automorphism groups of $V$, $Q$ and $S$. 
We begin by describing $\Aut(S)$.

\begin{lemma}\label{OutS} 
Let $A := \Aut_{\Aut(SU_4(p))}(S)$. Then the following hold:

\begin{enumerate}
\item The order of $\Aut(S)$ is $2p^a(p+1)(p-1)^2$ for some $a\in \mathbb{Z}_{\geq 0}$.
\item $A\cong \Inn(S)\rtimes(C_{p-1}\times (C_{p^2-1} \rtimes C_2))$.
\item Every subgroup of $\Aut(S)$ containing $\Inn(S)$ with $p'$-index is conjugate to a subgroup of $A$.
\item $\Aut_\FF(S)$ is $\Aut(S)$-conjugate to a subgroup of $A$ for any saturated $\FF$ on $S$.

\end{enumerate}
\end{lemma}

\begin{proof}Consider the chain $\mathcal{C}:\Phi(S)\trianglelefteq V\trianglelefteq S$ of characteristic subgroups of $S$. The stabilizer of $\mathcal{C}$ is a normal $p$-subgroup of $\Aut(S)$ by Lemma \ref{Oliver1.7}, and any other element of $\Aut(S)$ acts nontrivially on $\mathcal{C}$. We consider the action of $\Aut(S)$ on $S/V$.

 Lemma \ref{LUPropertiesOfS}~(6) implies that all $p+1$ elements of $\XX$
 are isomorphic with $M'\neq Z(M)$ and there is an element of order $p+1$ permuting them. This element acts transitively on 
 the $p+1$ nontrivial proper subgroups of $S/V\cong C_p^2$, hence its only overgroups in $GL_2(p)$ either contain $SL_2(p)$ or are contained in $C_{p^2-1}\rtimes C_2$, the normalizer in $GL_2(p)$ of a Singer cycle by \cite[II.7.3 and II.8.5]{Huppert}.
 There are no $p$-elements in $\Aut(S)/C_{\Aut(S)}(\mathcal{C})$, as one such would normalize some $M\in \XX$ and permute transitively the remaining 
 elements $N$ of $\XX$.  Hence the $p$-element would normalize $C_S(M')\in \XX\setminus \{M\}$, a contradiction.  In particular, $\Aut(S)/C_{\Aut(S)}(\mathcal{C})$ embeds into $GL_1(p)\times GL_2(p)$. 
Hence $|\Aut(S)|_{p'}$ divides  $2(p-1)^2(p+1)$. To obtain equality we observe that in $\Aut(SU_4(p))$ we have 
$|\Out_{\Aut(SU_4(p))}(S)| = 2(p-1)^2(p+1)$ 
by \cite[Table 8.10]{BHRD} and \cite[Table 2.1.C]{KLClassical}, hence $|\Aut(S)| = 2p^a(p+1)(p-1)^2$ as claimed and, as $\Inn(S)$ is characteristic in $\Aut(S)$, the isomorphism type of $A$ is $\Inn(S)\rtimes (C_{p-1}\times (C_{p^2-1} \rtimes C_2))$, proving parts (1) and (2).

Hence, $\Aut(S)$ is solvable and Hall's Theorem \cite[Theorem 6.4.1]{Gorenstein} implies that every subgroup of $\Aut(S)$ containing $\Inn(S)$ with $p'$-index is conjugate to a subgroup of $A$. In particular this is true for $\Aut_\FF(S)$, concluding the Lemma.
\end{proof}

We now determine uniqueness of two subgroups of $\Aut_{\Aut(SU_4(p))}(S)$, which will end up being $\Aut_\FF^0(S)$ for various values of $p$.

\begin{lemma}\label{UniquenessOfOutF0S} There is a unique subgroup $H$ of $\Out_{\Aut(SU_4(p))}(S)$ that is isomorphic to $C_{p-1}\times C_{(p^2-1)/2}$. Furthermore, $H$ has two subgroups isomorphic to $C_{(p-1)/2}\times C_{(p^2-1)/2}$, only one of which contains an element of order $(p^2-1)/2$ acting via an element of order $p-1$  on $S/Q$.
\end{lemma}

\begin{proof}
By Lemma \ref{OutS}~(2), we have $\Out_{\Aut(SU_4(p))}(S)\cong C_{p-1}\times D_{2(p^2-1)}$, so we fix generators $x,y,z$ such that 
$$\Out_{\Aut(SU_4(p))}(S) = \langle x\mid x^{p-1}\rangle \times \langle y,z \mid  y^{p^2-1}, z^2, y^z = y^{-1}\rangle.$$ 
As in Lemma \ref{OutS}, we see that $x$ acts on $S/Q$ and centralizes $Q/\Phi(S)$ whereas $\langle y, z\rangle$ centralizes $S/Q$ and acts on $Q/\Phi(S)\cong S/V$ as the normalizer in $GL_2(p)$ of a Singer cycle, which is a dihedral group by \cite[II.8.4]{Huppert}. We consider subgroups of index $4$, noting that $H = \langle x, y^2\rangle$ is one such. 

We also see that, since $p>3$, $z$ does not centralize $y^{(p^2-1)/4}$, hence $z$ is not contained in the subgroups in question.
There are three subgroups of index $2$ in $\langle x, y\rangle$, which contain $\langle x^2, y^2\rangle$: $H_1 = \langle x^2, y\rangle$, $H_2 = \langle x, y^2\rangle$ and $H_3 = \langle x^2, y^2, xy\rangle$. Now $H_1$ and $H_3$ contain an element of order $p^2-1$, hence are not of the required shape and thus $H = H_2$ is the unique subgroup of the given isomorphism type. 

We now consider subgroups $K_i$ of index $2$ in $H$, which must contain $\langle x^2, y^4 \rangle$.  Hence there are again 3 such: 
 $K_1 = \langle x^2, y^2\rangle$, $K_2 = \langle x^2, y^4, xy^2\rangle$ and $K_3 = \langle x, y^4\rangle$. Since the ones we are interested are $K_i\cong C_{p-1}\circ_{C_2} C_{(p^2-1)/2}$, such $K_i$ contain an element of order $(p^2-1)/2$; hence $K_i\cong C_{(p-1)/2}\times C_{(p^2-1)/2}$. However, $K_3$ has exponent $(p^2-1)/4$; hence it is not an option, but both $K_1$ and $K_2$ are isomorphic to $C_{(p-1)/2}\times C_{(p^2-1)/2}$, as required. 
If we further assume that we have an element which acts on $S/Q$ as an element of order $p-1$, we observe that in $K_1$ there is no such element. Thus $K_2$ is the only subgroup with the desired property, since the element $xy^2\in K_2$ acts on $S/\Phi(S)$ as desired, and the lemma is complete.
\end{proof}

Regarding $Q$, the automorphism group of an extraspecial $p$-group of exponent $p$ is well-known.

\begin{thm}[{\cite[Theorem 1]{Winter}}]
\label{AutomorphismOfExtraspecial}
Let $Q$ be an extraspecial group of order $p^{1+2n}$ and exponent $p$. Denote by $A = \Aut(Q)$, $B=C_A(Z(Q))$, and $C= C_B(Q/Z(Q))$. Then we have:
\begin{enumerate}
\item $ C=\Inn(Q)\cong Q/Z(Q)$ is elementary abelian of order $p^{2n}$;

\item $A = BT$ is the semidirect product of $B$ with a cyclic group $T$ of order $p-1$;

\item $B/C \cong Sp_{2n}(p)$. 
\end{enumerate}
\end{thm}

The following result will help us determine the action of $O^{p'}(\Aut_\FF(Q))$ on $Q/Z(Q)$.

\begin{prop}\label{pNot3V4dQuadratic} Assume $p > 3$ and let $G$ be a group with $S\in \Syl_p(G)$ of order $p$, $O_p(G) = 1$. Let $V$ be a $4$-dimensional faithful $GF(p)G$-module with $C_V(S) = [V,S]$ of dimension $2$. Then, $O^{p'}(G) \cong SL_2(p)$ and $V$ is a direct sum of two natural $SL_2(p)$-modules.
\end{prop}

\begin{proof}
Since $G$ acts faithfully on $V$, $G$ embeds into $\Aut(V)\cong GL_4(p)$.  Assume $G$ is a minimal counterexample to the lemma, that is, $S\leq G\leq GL_4(p)$ with $|G|$ minimal such that $O_p(G) = 1$ and if $S\leq L< G$ with $O_p(L) = 1$ then $O^{p'}(L)\cong SL_2(p)$ and $V$ is a direct sum of two natural $SL_2(p)$-modules for $L$.

Let $S = \langle s\rangle$. Since $C_V(S) = [V,S]$, we have $[V,s,s] = 0$ and, as $O_p(G) = 1$, $s\in G\setminus O_p(G)$.
Because $V$ is a faithful $G$-module,  \cite[Lemma 2.4]{ChermakQuadratic} yields  $G$ has a subgroup $H$ such that $H = \langle s^H\rangle \cong SL_2(p)$ and $V = [V,H]\oplus C_V(O^{p}(H))$, where $[V,H]$ is a direct sum of natural $SL_2(p)$-modules for $H$.
 Since $p>3$, $O^p(H) = H$, so $C_V(O^p(H)) =0$,  as otherwise  $\dim([V,s])=1$.  Hence, $V=[V,H]$ has dimension $4$.
 Furthermore, as an $H$-module, $V = [V,H]$ is a direct sum of two natural $SL_2(p)$-modules. In particular,  the central involution $t\in Z(H)$ negates $V$ and so  $t = -I_4\in Z(GL_4(p))$ and $t\in Z(G)$.
Therefore, $$C_V(G)\le C_V(H) \le C_V(t)=0.$$
 Let $T \in \Syl_p(H) \setminus\{S\}$, then $H= \langle S, T\rangle.$  Since $G$ is a counterexample to the lemma,  $H \ne O^{p'}(G)$ and so $H$ is not normal in $G$.

 Let $U \in \Syl_p(G) \setminus\{S\}$ and set $L = \langle S, U\rangle$. Since $S\ne U$, $O_p(L) = 1$.
Assume that $C_V(S) \cap C_V(U) \neq 0$. Then, $$C_V(L) = C_V(S) \cap C_V(U)>0=C_V(G)$$ and so $L <G$.
Since  $L$ satisfies the hypothesis of the lemma, induction implies that $L=O^{p'}(L)\cong SL_2(p)$ and $V|_L$ is a direct sum of two natural $SL_2(p)$-modules for $L$. But this means $C_V(L) =0$, a contradiction.

 Therefore, if $R, U \in \Syl_p(G) $ with $R \ne U$, then \begin{equation}\label{Centralizers}\tag{\textdagger}C_V(R) \cap C_V(U) = 0.\end{equation} In particular, \begin{equation*}
 |\bigcup_{R \in \Syl_p(G)}C_V(R)|= |\Syl_p(G)|(p^2-1) +1 \le p^4.\end{equation*}
Hence $\left | \mathrm{Syl}_p(G)\right|\leq p^2+1.$

 We  investigate $N_G(S)$.  Since $\Aut(S)$ is abelian, we have $N_G(S)' \le C_G(S).$  Furthermore, $N_G(S)$ acts on $C_V(S)$ and so $N_G(S)/C_{N_G(S)}(C_V(S))$ is isomorphic to a subgroup of $\Aut(C_V(S))$, which is isomorphic to $GL_2(p)$.  By Lemma \ref{AuxiliaryMeier1}, the eigenvalues of elements of $N_H(S)$ on $C_V(S)$ are equal, and so $N_H(S)C_{N_G(S)}(C_V(S))/C_{N_G(S)}(C_V(S))$ acts on $C_V(S)$ as scalars. 
 %(the eigenvalues of elements of $N_H(S)$ on $C_V(S)$ are equal, see \cite[Lemma 1.48]{RMMThesis}).  
 It follows that $$[N_H(S),N_G(S)] \le  C_{N_G(S)}(C_V(S))\cap N_G(S)'\le  C_{N_G(S)}(C_V(S))\cap C_G(S).$$

 Assume that $x\in [N_H(S),N_G(S)]$ has $p'$-order. Then $[V,S,x]= [[V,S],x]=0$ and so the Three Subgroups Lemma implies $[V,x,S]=0$.  Hence $[V,x, x]=0$ and so $x$ centralizes $V$ by coprime action.  It follows by \cite[Theorem 5.3.6]{Gorenstein} that $x=1$. Hence  $[N_H(S),N_G(S)]$ is a $p$-group.  As $[N_H(S),N_H(S)]=S$ and $S \in \Syl_p(G)$, we conclude that $$[N_H(S),N_G(S)]= S\le N_H(S)$$ and so $N_H(S)$ is normal in $N_G(S)$.

Assume that no two distinct conjugates of $H$ contain $S$.  Then, by Sylow's Theorem, for $k, \ell \in G$ with $H^k \ne H^\ell$, $|H^k \cap H^\ell|$ is coprime to $p$.   Suppose $K$ is a conjugate of $H$ with $H \ne K$.
 Then, as $S \le H$ and $S \in \Syl_p(G)$, $S$ does not normalize $K$, as otherwise $p^2\mid |SK|$. Hence $|\{K^s\mid s \in S\}|=p$ and, for all $s\in S$, $p$ does not divide $|K^s \cap K|$. Thus
\begin{eqnarray*}|\Syl_p(G) | &\ge& |\Syl_p(H) | + \sum _{s \in S}|\Syl_p(K^s)|=p+1+p(p+1)\\& = & p^2+2p +1 > p^2+1 \ge |\Syl_p(G)|,\end{eqnarray*} a contradiction. Thus, there exist two distinct conjugates of $H$ containing $S$.

 Let $K$ be one such, that is a conjugate of $H$ with $K\neq H$ and $S\leq H\cap K$. Then $K= H^g$ for some $g \in G$ and so $S, S^g \le K$. By Sylow's Theorem, there exists $k \in K$ such that $S^{gk}= S$. Now $H^{gk}= K^k=K$.  Hence, we may assume that $g \in N_G(S)$. In particular, as $N_H(S) $ is normal in $N_G(S)$, $N_H(S)=N_H(S)^g\le H^g =K$.  Thus, $N_H(S)= N_K(S)$.  Let $X$ be a complement to $S$ in $N_H(S)$ which normalizes $T$. Then, $X$ is cyclic of order $p-1$, $X\le N_K(S)$ and $X$ normalizes some  $U\in \Syl_p(K)\setminus\{S\}$.  Let $x$ be a generator of $X$ and note that, by the assumption $p>3$, $x$ is not an involution.
Therefore, Lemma \ref{AuxiliaryMeier2} 
%a straightforward calculation (see \cite[Lemma 1.49]{RMMThesis}) 
yields that $x$ has exactly two eigenvalues $\lambda$ and $\lambda^{-1}$ on $V$ and the corresponding eigenspaces are $C_V(S)$ and $C_V(T)$. Since $C_V(S) \cap C_V(U) =0$ and $X$ acts on $C_V(U)$, there is an eigenvector for $x$ in $C_V(U)$ which is not in $C_V(S)$. It follows that $C_V(U) \cap C_V(T) \ne 0$. By (\ref{Centralizers}), we conclude that $T= U$.
 But then $K=\langle S, U\rangle=\langle S,T\rangle =H,$
a contradiction. This contradiction proves that $G$ is not a counterexample and proves the lemma.
\end{proof}

Finally, we turn to $V$, which is elementary abelian, so we have $\Aut(V)\cong GL_4(p)$. The following result gives the information which we use to determine  $O^{p'}(\Aut_\FF(V))$.

\begin{prop}\label{speinGL4}
Suppose $p\geq 5$, $V\cong C_p^4$ and  $G\leq \Aut(V)\cong GL_4(p)$ has a strongly $p$-embedded subgroup and $G$ has $p$-rank at least $2$. Then, $O^{p'}(G)$ is isomorphic to either $SL_2(p^2)$ or $PSL_2(p^2)$. Furthermore, if a Sylow $p$-subgroup of $G$ fixes a $1$-subspace of $V$, then $O^{p'}(G)\cong PSL_2(p^2)$ and $V$ is the natural $\Omega_4^-(p)$-module.
\end{prop}

\begin{proof}Under the assumptions above, Theorem \ref{almostsimplespe} implies that $K := O^{p'}(G/O_{p'}(G))$ is almost simple, and we
can use \cite[Theorem 7.6.1]{GLS3} to obtain a list of candidates for $K$. Let $T\in \Syl_p(K)$. We first rule out all candidates except $PSL_2(p^2)$, then show that $O_{p'}(G)$ centralizes $O^{p'}(G)$, so that the result follows.

\begin{claim}$K\cong PSL_2(p^2)$.
\end{claim}
\begin{claimproof}[Proof of claim]
Note that a Sylow $p$-subgroup $U$ of $GL_4(p)$ has order $p^6$, nilpotency class $3$, and $K$ must be isomorphic to a section of $GL_4(p)$, in particular, $|K|$ must divide $$|GL_4(p)| = p^6(p^4-1)(p^3-1)(p^2-1)(p-1).$$ We obtain the candidates for $K$ from \cite[Theorem 7.6.1]{GLS3} and the structure of $T$ is an easy consequence, see e.g. \cite[Corollary 1.60]{RMMThesis}. We prove that they do not embed using divisibility arguments, which also cover the case of their covers.

If $K/Z(K)\cong PSU_3(p^n)$, with $p^n\geq 3$, then $T$ has order $p^{3n}$, so we must have $n\leq 2$. If $n=2$ then $T$ has order $p^6$, so we need $U\cong T$, but $T$ has nilpotency class $2$, so it cannot happen. Finally, if $n=1$, we have $p^3+1\mid|PSU_3(p)|$, and as $p^6-1 = (p^3+1)(p^3-1)$, Zsigmondy's Theorem  \cite{Zsigmondy}
gives a contradiction.

If $K/Z(K)\cong Sz(2^{2n+1})$ then $|T| = (2^{2n+1})^2\leq 2^6$, so we must have $n=1$, but then $T\cong U$. However, $T$ has nilpotency class $2$ while $U$ has nilpotency class $3$, so this case does not happen either.

If $K/Z(K)\cong {^2G_2(3^{2n+1})}$ then $|T| = (3^{2n+1})^3 > 3^6$, a contradiction.

If $K/Z(K)\cong A_{2p}$, as $p\geq 5$, then we must have $p\leq 17$, as otherwise we observe that $|A_{2p}| = (2p)!/2\geq p^{17} > p^{16} \geq (p^4-1)(p^3-1)(p^2-1)(p-1)p^6 = |GL_4(p)|$, so it cannot embed. For the remaining primes $p=5,7,11,13$, we have, respectively, the primes $q = 7,13,17,23$ such that $q\mid |A_{2p}|$ but $q\nmid |GL_4(p)|$, so this case does not happen.

 For the remaining cases other than $PSL_2(p^n)$ there is always a prime dividing $|K|$ that does not divide $|GL_4(p)|$. 
 If $p=5$, none of $11,13,41$ divide $|GL_4(5)|$ but $11$ divides $|McL|$ and $|Fi_{22}|$, $13 \mid |^2F_4(2)'|$, and $41\mid |\Aut(Sz(32))|$. And for $p=11$, we have that $43$ does not divide the order of the sporadic Janko group $J_4$, but $43\nmid |GL_4(11)|$.

Finally, if $K\cong PSL_2(p^n)$ then $|K|= (p^{2n}-1)p^n/2 $ by \cite[Theorem 2.8.1]{Gorenstein}. 
Then, if $n\geq 3$, we have again by Zsigmondy's Theorem \cite{Zsigmondy}
a prime $q$ such that $q\mid p^{2n}-1$ and $q\nmid p^k-1$ for any $k<2n$ unless $p=2$ and $n=3$. The Sylow $p$-subgroups of $PSL_2(p)$ are cyclic, so we cannot have $n=1$. 
Therefore, we must have $K\cong PSL_2(p^2)$, which embeds into $GL_4(p)$ as $PSL_2(p^2)\cong \Omega^-_4(p)$ by \cite[Proposition 2.9.1~(5)]{KLClassical}.
\end{claimproof}

\begin{claim}
$O_{p'}(G)\leq Z(O^{p'}(G))$.
\end{claim}

\begin{claimproof}[Proof of claim]We have $K=O^{p'}(G/O_{p'}(G))\cong PSL_2(p^2)$. 
Let $R := O_{p'}(G)$ and $T\in \Syl_p(G)$. Then, $T\cong C_p\times C_p$, so we pick $x\in T$ of order $p$ and
consider {$H := R\langle x\rangle\leq G$}. $H$ is $p$-solvable as $R\trianglelefteq H$ is a $p'$-group and $H/R\cong \langle x\rangle$ is a $p$-group.
Consider the action of $x$ on the natural $GL_4(p)$-module. The Jordan form of $x$ has largest Jordan block of size at most $4$, so that its minimal polynomial is $(X-1)^r$ for some $r\leq 4$.

Then, if $O_p(H) = 1$, we have by the Hall-Higman Theorem \cite[Theorem 11.1.1]{Gorenstein} that $p-1\leq r\leq p$.
 This means that if $p\geq 7$ or $p=5$ and $r\leq 3$ then $O_p(H)\neq 1$. Thus in these cases $H = R\times \langle x\rangle$, and $R$ centralizes $x$. Since we can do this for any subgroup of any Sylow $p$-subgroup of $G$, we conclude that $R$ acts trivially on $O^{p'}(G)$ and therefore $R \leq Z(O^{p'}(G))$ and $O^{p'}(G)$ is a central extension of $PSL_2(p^2)$.
 The Schur multiplier of $PSL_2(p^2)$ has order $2$  if $p\geq  5$ and its universal covering group is $SL_2(p^2)$ by 
\cite[V.25.7 Satz]{Huppert}, so in this case the result follows.

The remaining case
 is $p=5$ and $x\in T$ has Jordan form $J_4$.

Let $C := C_{GL_4(5)}(x)$, then $T\leq C$ and $C$ has shape $C_5^3 : C_4$ by \cite[Theorem 7.1]{LiebeckSeitzUnipotent}. We claim there is no subgroup of the centralizer of order $5^2$ with only elements with Jordan form $J_4$. This is because a Sylow $5$-subgroup of $C$ is generated by matrices 

\begin{center} $x = \left(\begin{array}{cccc} 1 & 0 & 0 & 0 \\
1 & 1 & 0 & 0 \\
0 & 1 & 1 & 0 \\
0 & 0 & 1 & 1 \end{array}\right)$, 
~ $x_2= \left(\begin{array}{cccc} 1 & 0 & 0 & 0 \\
0 & 1 & 0 & 0 \\
1 & 0 & 1 & 0 \\
0 & 1 & 0 & 1\end{array}\right)$, 
~ $x_3= \left(\begin{array}{cccc} 1 & 0 & 0 & 0 \\
0 & 1 & 0 & 0 \\
0 & 0 & 1 & 0 \\
1 & 0 & 0 & 1 \end{array}\right)$, 
\end{center}
so that the subgroup generated by $x_2$ and $x_3$ contains no element with Jordan form $J_4$. Any subgroup of order $p^2$ in $C_{GL_4(5)}(x)$ must intersect this subgroup nontrivially, hence it must contain some nonidentity element $y$ with Jordan form distinct from $J_4$.
Then, as before, $R$ must centralize $y$ by Hall-Higman. Note that the subgroup $\langle y^G\rangle =O^{p'}(G)$, since $O^{p'}(G/O_{p'}(G))$ is almost simple,
so that $R$ centralizes $O^{p'}(G)$ in this case as well. Thus, in every case, $O^{p'}(G)$ is a central extension of $PSL_2(p^2)$.
\end{claimproof}

The first part of the proposition now follows, and both cases arise.

Assume further that $\dim C_V(S)=1$, where $S\in \Syl_p(G)$, then $V$ is not the natural $SL_2(p^2)$-module and the  description of the $GF(p)$-modules for $SL_2(p^2)$ in \cite[Section 30]{BrauerNesbitt} implies that $V$ is the natural $\Omega_4^-(p)$-module.
\end{proof}

We now describe the relevant parts of the structure of the natural $\Omega_4^-(p)$-module.

\begin{lemma}\label{Omega4-Calculation} Let $p$ be odd. Suppose $G\cong \Omega_4^-(p)\leq GL_4(p)$ acts on the natural $\Omega_4^-(p)$-module $V$,  let $R\in \Syl_p(G)$ and $K=\langle t\rangle$ be a complement to $R$ in $N_G(R)$. Then $|R| = p^2$,
\begin{enumerate}
\item $R$ preserves exactly $p$ non-degenerate quadratic forms on $V$ up to scalars.
\item  $V = V_1\oplus V_2 \oplus V_3$ as a $K$-module with $V_i$ irreducible, $\dim(V_1) = \dim(V_3) = 1$ and $\dim(V_2) = 2$. $ V_3 = C_V(R)$ and $[V,R] = V_2\oplus V_3$.
The element $t$ has order $(p^2-1)/2$ and acts as an element of order $p-1$ on $V_1$ and $V_3$, and as an element of order $p+1$ on $V_2$. 

\item $C_K(C_V(R)) = \langle t^{p-1}\rangle\cong C_{(p+1)/2}$.

\item Let $i = t^{(p^2-1)/4}$ be the unique involution in $K$. 
If $4\mid p+1$ then $i$ centralizes $V_1$ and $V_3$ and inverts $V_2$, whereas if $4\mid p-1$ then $i$ centralizes $V_2$ and inverts $V_1$ and $V_3$. In either case, $i$ inverts $R$.

\item If $p\geq 5$ then there is a unique non-degenerate quadratic form up to a scalar which is preserved by both $R$ and $t$, that is by $N_G(R)$.
\end{enumerate}

\end{lemma}
\begin{proof} We prove the result by explicit calculation. 
%, which can be found in \cite[Lemma 5.7]{RMMThesis}
%\end{proof}
%
%-------------
%
%
%\section{Natural $\Omega_4^-(p)$-module calculations}
%
%
%
%\begin{lemma}\label{Omega4-Calculation} Suppose $G\cong \Omega_4^-(p)\leq GL_4(p)$ acts on the natural $\Omega_4^-(p)$-module $V$,  let $R\in \Syl_p(G)$ %consists of lower triangular matrices, 
%and $K=\langle t\rangle$ be a complement to $R$ in $N_G(R)$. Then $|R| = p^2$,
%\begin{enumerate}
%\item $R$ preserves exactly $p$ non-degenerate quadratic forms on $V$ up to scalars.
%\item  $V = V_1\oplus V_2 \oplus V_3$ as a $K$-module with $V_i$ irreducible, $\dim(V_1) = \dim(V_3) = 1$ and $\dim(V_2) = 2$. $ V_3 = C_V(R)$ and $[V,R] = V_2\oplus V_3$.
%The element $t$ has order $(p^2-1)/2$ and acts as an element of order $p-1$ on $V_1$ and $V_3$, and as an element of order $p+1$ on $V_2$. 
%
%% irreducibly on $[V,R]/C_V(R)$ and again via an element of order $p-1$ on $V/[V,R]$. 
%
%%, $V_1 = C_V(R)$, and $[V,R] = C_V(R)\oplus V_2$.
%\item $C_K(C_V(R)) = \langle t^{p-1}\rangle\cong C_{(p+1)/2}$.
%
%\item Let $i = t^{(p^2-1)/4}$ be the unique involution in $K$. 
%%The action of the unique involution in $i := t^{(p^2-1)/4}\in K$ on $V$ depends on the value of $p\pmod 4$. 
%If $4\mid p+1$ then $i$ centralises $V_1$ and $V_3$ and inverts $V_2$, whereas if $4\mid p-1$ then $i$ centralises $V_2$ and inverts $V_1$ and $V_3$. In either case $i$ inverts $R$.
%
%
%\item If $p\geq 5$ then there is a unique non-degenerate quadratic form up to a scalar which is preserved by both $R$ and $t$.
%\end{enumerate}
%
%\end{lemma}
%
%\begin{proof}
Since $G\cong \Omega_4^-(p)$, it leaves invariant a quadratic form with matrix $F$. Then, as $G$ has type $-$, $F$ has Witt index $1$ and \cite[21.2]{aschbacher} implies that $V = O D$ where $O$ is a 2-dimensional definite orthogonal space and $D$ a hyperbolic plane. Therefore, there is a basis $\{v_1, v_2, v_3, v_4\}$ of $V$ such that \cite[19.2 and 21.1]{aschbacher} show that $F =  \small\left(\begin{array}{cccc} 0 & 0 & 0 & 1 \\
0 & 1 & 0 & 0 \\
0 & 0 & \alpha & 0 \\
1 & 0 & 0 & 0 \end{array}\right)$\normalsize, where $-\alpha$ generates $GF(p)$, so that every $g\in G$ satisfies $gFg^T = F$ by \cite[19.7]{aschbacher}. Fix this basis. %There is a unique form up to "iso???" so that there is a unique conjugacy class in $GL_4(p)$ of such subgroups. Hence we may assume that a Sylow $p$-subgroup $R$ of $G$ consists of lower triangular matrices. 
Recall that, as $\Omega_4^-(p)\cong PSL_2(p^2)$ by 
%Proposition \ref{KL291}~(5)
\cite[Proposition 2.9.1~(5)]{KLClassical}, we have $|R|= p^2$ and $N_G(R) =R\rtimes K\cong  C_p^2 \rtimes C_{(p^2-1)/2}$, and we can find $R$ as lower triangular matrices.
We now calculate which lower triangular matrices preserve the quadratic form given by $F$. Let $a,b,c,d,e,f\in GF(p)$. Then \noindent
\begin{align*} F =  \small \left(\begin{array}{cccc} 0 & 0 & 0 & 1 \\
0 & 1 & 0 & 0 \\
0 & 0 & \alpha & 0 \\
1 & 0 & 0 & 0 \end{array}\right) 
=\left(\begin{array}{cccc} 1 & 0 & 0 & 0 \\
a & 1 & 0 & 0 \\
b & c & 1 & 0 \\
d & e & f & 1 \end{array}\right)  
\left(\begin{array}{cccc} 0 & 0 & 0 & 1 \\
0 & 1 & 0 & 0 \\
0 & 0 & \alpha & 0 \\
1 & 0 & 0 & 0 \end{array}\right)
\left(\begin{array}{cccc} 1 & a & b & d \\
0 & 1 & c & e \\
0 & 0 & 1 & f \\
0 & 0 & 0 & 1 \end{array}\right) = \\
\left(\begin{array}{cccc} 0 & 0 & 0 & 1 \\
0 & 1 & 0 & a \\
0 & c & \alpha & b \\
1 & e & \alpha f & d \end{array}\right)
\left(\begin{array}{cccc} 1 & a & b & d \\
0 & 1 & c & e \\
0 & 0 & 1 & f \\
0 & 0 & 0 & 1 \end{array}\right) =
\left(\begin{array}{cccc} 0 & 0 & 0 & 1 \\
0 & 1 & c & e+a \\
0 & c & c^2+\alpha & ec+\alpha f + b \\
1 & a+e & b+ce+\alpha f & d+e^2+\alpha f^2+d \end{array}\right)
\end{align*}
\normalsize
therefore any such matrix satisfies $c=0$, $e=-a$, $b = -\alpha f$, and $d=-(a^2+f^2)/2$, so that we have  generators  $r_1 =  \small\left(\begin{array}{cccc} 1 & 0 & 0 & 0 \\
1 & 1 & 0 & 0 \\
0 & 0 & 1 & 0 \\
-1/2 & -1 & 0 & 1 \end{array}\right)$\normalsize and $r_2 = 
\small\left(\begin{array}{cccc} 1 & 0 & 0 & 0 \\
0 & 1 & 0 & 0 \\
-\alpha & 0 & 1 & 0 \\
-\alpha/2 & 0 & 1 & 1 \end{array}\right)$\normalsize of $R$. Hence, $R = \langle r_1, r_2\rangle$ is the group of lower triangular matrices preserving $F$, and by comparing orders we see that $R\in \Syl_p(G)$. 
Note that we can see that, with respect to this basis, $ C_V(R) = \langle v_4 \rangle$ and $[V, R] = \langle v_2, v_3, v_4\rangle$. We denote $V_1 = \langle v_1\rangle$, $V_2 = \langle v_2, v_3\rangle$ and $V_3 = \langle v_4\rangle = C_V(R)$. In particular, $[V, R] = V_2\oplus V_3$.

We now consider how many non-degenerate quadratic forms $R$ leaves invariant. Suppose there is a quadratic form $L$ which is invariant under $R$, then we have $L = r_i L r_i^T$ for $i=1,2$, hence

\noindent
\begin{align*} L=
\small\left(\begin{array}{cccc} 
a & b & c & d \\
b & e & f & g \\
c & f & h & i \\
d & g & i & j \end{array}\right)
%=
%\left(\begin{array}{cccc} 1 & 0 & 0 & 0 \\
%1 & 1 & 0 & 0 \\
%0 & 0 & 1 & 0 \\
%-1/2 & -1 & 0 & 1 \end{array}\right)  
%\left(\begin{array}{cccc} a & b & c & d \\
%b & e & f & g \\
%c & f & h & i \\
%d & g & i & j \end{array}\right)
%\left(\begin{array}{cccc} 1 & 1 & 0 & -1/2 \\
%0 & 1 & 0 & -1 \\
%0 & 0 & 1 & 0 \\
%0 & 0 & 0 & 1 \end{array}\right) &\\
%=
%\left(\begin{array}{cccc} a & b & c & d \\
%a+b & b+e & c+f & d+g \\
%c & f & h & i \\
%-a/2-b+d & -b/2-e+g & -c/2-f+i & -d/2-g+j \end{array}\right)
%\left(\begin{array}{cccc} 1 & 1 & 0 & -1/2 \\
%0 & 1 & 0 & -1 \\
%0 & 0 & 1 & 0 \\
%0 & 0 & 0 & 1 \end{array}\right) & \\
=
\left(\begin{array}{cccc} a & a+b & c & -a/2-b+d \\
a+b & a+2b+e & c+f & x_1 \\
c & c+f & h & -c/2-f+i \\
-a/2-b+d & x_1 & -c/2-f+i & x_2 \end{array}\right)
\normalsize= r_1 L r_1^T
%=
% \left(\begin{array}{cccc} a & b & c & d \\
%b & e & f & g \\
%c & f & h & i \\
%d & g & i & j \end{array}\right)&
\end{align*}

\noindent where $x_1 = -a/2-b+d-b/2-e+g$, $x_2 = -d/2+e-g-d/2-g+j$. We see that $a+b=b$, so $a=0$, then $a+2b+e = 2b+e = e$, so $b=0$. Further, since $c+f=f$, we have $c=0$ too, and as $-c/2-f+i = i$, we see that $f=0$ as well. As $x_1 = d-e+g= g$ we deduce $d=e$, and as $x_2 = -d+e-2g+j=j$, we see that $g=0$, so that 
 $L =  
 \small\left(\begin{array}{cccc} 
0 & 0 & 0 & d \\
0 & d & 0 & 0 \\
0 & 0 & h & i \\
d & 0 & i & j \end{array}\right)$\normalsize.
To determine $h,i,j$ we compute %perform a very similar computation 
$L = r_2 L r_2^T =  \small\left(\begin{array}{cccc} 
0 & 0 & 0 & d \\
0 & d & 0 & 0 \\
0 & 0 & h & h-\alpha d + i \\
d & 0 & -\alpha d+h+i & -\alpha d/2+h+i-\alpha d/2+i+j \end{array}\right)$\normalsize,
which shows that
$h=\alpha d$, $i=0$, and puts no restrictions on $j$. Hence, the non-degenerate quadratic forms fixed by $R$ have matrix
$L =  \small\left(\begin{array}{cccc} 
0 & 0 & 0 & d \\
0 & d & 0 & 0 \\
0 & 0 & \alpha d & 0 \\
d & 0 & 0 & j \end{array}\right)$\normalsize, and up to scalars there are $p$ such forms %which differ by an element of the centre of a Sylow $p$-subgroup of $GL_4(p)$, 
Part (1) holds.

We now consider $N_G(R) =  R\rtimes K$ where $K$ is cyclic of order $(p^2-1)/2$, and we claim that $K = \langle t\rangle$ where
$t = \small \left(\begin{array}{cccc} 
\lambda & 0 & 0 & 0 \\
0 & a & b & 0 \\
0 & c & e & 0 \\
0 & 0 & 0 & \gamma \end{array}\right)$\normalsize
satisfying $\lambda \gamma = 1$ with $\lambda$ a primitive element of $GF(p)$. We note that as $t\in G$, $\det(t) = 1$, so that the submatrix $M = \left(\begin{smallmatrix}a & b\\c&e\end{smallmatrix}\right)$ preserves the quadratic form $N = \left(\begin{smallmatrix}1 & 0\\0&\alpha\end{smallmatrix}\right)$ of type $-$ and has determinant $1$, hence embeds into $SO_2^-(p)\cong C_{p+1}$ by 
%Proposition \ref{KL291}~(3)
\cite[Proposition 2.9.1~(3)]{KLClassical}, and  can be chosen to have order $p+1$. In order to preserve the form $N$, it must  satisfy $a^2+b^2\alpha = 1$, $ca+be\alpha = 0$ and $c^2+e^2\alpha = \alpha$. 
Hence, we have  $t^{(p+1)(p-1)/\gcd(p+1, p-1)} = 1$ and $t$ has order $(p^2-1)/2$.

We see $t\in G$ when the above holds as it preserves $F$:  
$$tFt^T =  \small\left(\begin{array}{cccc} 
0 & 0 & 0 & \lambda\gamma \\
0 & a^2+b^2\alpha & ca+be\alpha & 0 \\
0 & ca+be\alpha & c^2+e^2\alpha & 0 \\
\lambda\gamma & 0 & 0 & 0 \end{array}\right)\normalsize = F.$$

We see that $t$ normalises $R$ since % as follows:
\noindent\begin{align*}t^{-1}r_1 t =  
%\left(\begin{array}{cccc} 
%\lambda^{-1} & 0 & 0 & 0 \\
%0 & e & -b & 0 \\
%0 & -c & a & 0 \\
%0 & 0 & 0 & \gamma^{-1} \end{array}\right)
% \left(\begin{array}{cccc} 1 & 0 & 0 & 0 \\
%1 & 1 & 0 & 0 \\
%0 & 0 & 1 & 0 \\
%-1/2 & -1 & 0 & 1 \end{array}\right)
% \left(\begin{array}{cccc} 
%\lambda & 0 & 0 & 0 \\
%0 & a & b & 0 \\
%0 & c & e & 0 \\
%0 & 0 & 0 & \gamma \end{array}\right) = \\
% \left(\begin{array}{cccc} \lambda^{-1} & 0 & 0 & 0 \\
%e & e & -b & 0 \\
%-c & -c & a & 0 \\
%-\gamma^{-1}/2 & -\gamma^{-1} & 0 & \gamma^{-1} \end{array}\right)
% \left(\begin{array}{cccc} 
%\lambda & 0 & 0 & 0 \\
%0 & a & b & 0 \\
%0 & c & e & 0 \\
%0 & 0 & 0 & \gamma \end{array}\right)
%= 
\small\left(\begin{array}{cccc} 
1 & 0 & 0 & 0 \\
\lambda e & 1 & 0 & 0 \\
-\lambda c & 0 & 1 & 0 \\
-\lambda\gamma^{-1}/2 & -\gamma^{-1} a & -\gamma^{-1} b & 1 \end{array}\right)
\end{align*}\normalsize
%where $D = ae-bc$ and 
 is a lower triangular matrix in $G$ as $t,r_1\in G$, and therefore $t^{-1}r_1t\in R$ as claimed.
A similar calculation works for $r_2$. %(Also note $D = ae-bc\in \pm 1$ as in $O_2^-(p)$)
%QUESTION: Should I try to get the explicit form in the previous way of this matrix? 
Note that we get $ e = a$  and $c =-\alpha b$.

Now, %as $V_1 = C_V(R) = \langle v_4Ê\rangle$, 
$t\in N_G(R)$ acts on $V_1 = \langle v_1\rangle$ and on $V_3 = C_V(R)=\langle v_4\rangle$ as an element of order $p-1$, and as an element of order $p+1$ on $V_2$ of order $p^2$, hence irreducibly. We have shown that parts (2) and (3) hold.

We can further see that the action of the involution  $i := t^{(p^2-1)/4}$ on $V$ depends on the value of $p\pmod4$. If $4\mid p+1$ then $p-1\mid (p^2-1)/4$ and $p+1\nmid (p^2-1)/4$, hence $t^{(p^2-1)/4}$ centralises $V_1$ and $V_3$ while inverting $V_2$, whereas if $4\mid p-1$ the divisibility conditions are swapped, hence $t^{(p^2-1)/4}$ centralises $V_2$ while inverting $V_1$ and $V_3$.
Hence \small $i = \left(\begin{array}{cccc} 
\epsilon & 0 & 0 & 0 \\
0 & -\epsilon & 0 & 0 \\
0 & 0 & -\epsilon  & 0 \\
0 & 0 & 0 & \epsilon \end{array}\right)$\normalsize, where $\epsilon = 1$ if $4\mid p+1$ and $\epsilon =-1$ if $4\mid p-1$.  In either case we have $ir_k i = r_k^{-1}$ for $k=1,2$, which establishes part (4).%, and part (4) follows. %Note that in either case $-I_4\notin N_G(T)$.

We now consider which quadratic forms are preserved simultaneously by $t$ and $R$. In part (1), we showed that any non-degenerate quadratic form preserved by $R$ is of the form
$L =  \small\left(\begin{array}{cccc} 
0 & 0 & 0 & d \\
0 & d & 0 & 0 \\
0 & 0 & \alpha d & 0 \\
d & 0 & 0 & j \end{array}\right)$\normalsize, 
and we have
$t =  \small\left(\begin{array}{cccc} 
\lambda & 0 & 0 & 0 \\
0 & a & b & 0 \\
0 & c & e & 0 \\
0 & 0 & 0 & \gamma \end{array}\right)$\normalsize. 

Then,
% We calculate 
 %
% \small\noindent
$L = tLt^T =   
%\left(\begin{array}{cccc} 
%\lambda & 0 & 0 & 0 \\
%0 & a & b & 0 \\
%0 & c & e & 0 \\
%0 & 0 & 0 & \gamma \end{array}\right)
%\left(\begin{array}{cccc} 
%0 & 0 & 0 & d \\
%0 & d & 0 & 0 \\
%0 & 0 & \alpha d & 0 \\
%d & 0 & 0 & j \end{array}\right)
%\left(\begin{array}{cccc} 
%\lambda & 0 & 0 & 0 \\
%0 & a & c & 0 \\
%0 & b & e & 0 \\
%0 & 0 & 0 & \gamma \end{array}\right) = \\
%\left(\begin{array}{cccc} 
%0 & 0 & 0 & \lambda d \\
%0 & ad & \alpha bd & 0 \\
%0 & cd & \alpha de & 0 \\
%d\gamma & 0 & 0 & \gamma j \end{array}\right)
%\left(\begin{array}{cccc} 
%\lambda & 0 & 0 & 0 \\
%0 & a & c& 0 \\
%0 & b & e & 0 \\
%0 & 0 & 0 & \gamma \end{array}\right) 
%=
\small\left(\begin{array}{cccc} 
0 & 0 & 0 & \gamma\lambda d \\
0 & a^2d+\alpha b^2d & acd+\alpha bde & 0 \\
0 & acd+\alpha bde & c^2d+\alpha de^2 & 0 \\
\lambda\gamma d & 0 & 0 & \gamma^2 j \end{array}\right) 
$\normalsize.

\noindent The middle $2\times 2$ submatrix is $I_2$ by the constraints $t$ satisfies, and we  obtain that $j = \gamma^2 j$, which must hold for each nonzero $\gamma\in GF(p)$ so, whenever $p > 3$, we must have $j= 0$. 

Therefore, \small$K =  \left(\begin{array}{cccc} 
0 & 0 & 0 & d \\
0 & d & 0 & 0 \\
0 & 0 & \alpha d & 0 \\
d & 0 & 0 & 0 \end{array}\right)$\normalsize is the unique non-degenerate quadratic form that is preserved by $N_G(R)$, and $d$ can be chosen to be $1$, hence in these circumstances claim (5) holds.
\end{proof}

\section{\texorpdfstring{$\FF$}{\textit{F}}-essential subgroups}\label{sec:essential}

Throughout this section we assume that $p\geq 5$,  $S$ is a Sylow $p$-subgroup of $SU_4(p)$, and $\FF$ is a saturated fusion system on $S$, 
 and we take the notation previously established.
In this section we determine the subgroups of $S$ which can be $\FF$-essential, and prove the following:

\begin{thm}\label{essentials}\begin{enumerate}
\item The only potential $\FF$-essential subgroups of $S$ are $Q$ and $V$.  
\item If $Q$ is $\FF$-essential then $Q/Z$ as a $O^{p'}(\Out_\FF(Q))$-module is a direct sum of two natural $SL_2(p)$-modules. 
\item If $V$ is $\FF$-essential then $V$ is a natural $\Omega_4^-(p)$-module for $O^{p'}(\Aut_\FF(V))\cong PSL_2(p^2)$.
\item ${O_p(\FF) = 1}$ if and only if  both $Q$ and $V$ are $\FF$-essential.
\end{enumerate}
\end{thm}

We first consider which nonabelian subgroups of $S$ of order $p^4$ can be $\FF$-essential.

\begin{lemma}\label{EssentialCandidatesOrderp4}
Suppose $H\leq S$ is nonabelian of order $p^4$ and $\FF$-essential.{~}Then $H\cong p^{1+2}_+\times C_p$.
\end{lemma}

\begin{proof}
Let $H\leq S$ of order $p^4$ be nonabelian and $\FF$-essential, then $O_p(\Out_\FF(H))=1$.  If $H$ has maximal nilpotency class then
%
%---
%
%$\Phi(H)$ has order $p^2$ and is abelian, with 
%
$|\Phi(H)| = p^2$ and $|Z(H)| = p$, so 
$Z(H)<\Phi(H)$, so the characteristic subgroup $\gamma_1(H) = C_H(\Phi(H))$ (described just before \cite[Lemma 2.5]{Blackburn}) has index $p$ in $H$. Thus,
 there is a characteristic  chain of subgroups $\Phi(H)\trianglelefteq C_H(\Phi(H))\trianglelefteq H$, each of index $p$ in the next, that are all characteristic in $H$.
 %
 %CHAIN
 %
 %, and which are centralized by $\Out_S(H)\cong N_S(H)/H$, which is a Sylow $p$-subgroup of $\Out_\FF(H)$. 
 Since the quotients $H/C_H(\Phi(H))$ and $C_H(\Phi(H))/\Phi(H)$ have order $p$, $\Out_S(H)$ centralizes them, 
 hence $\Out_S(H)$ is a normal $p$-subgroup of $\Out_\FF(H)$ by Lemma \ref{Oliver1.7} and thus $O_p(\Out_\FF(H))\neq 1$, contrary to $O_p(\Out_\FF(H)) = 1$.

Thus, $H$ has nilpotency class $2$ and $H'\leq Z(H)$. If $Z(H) = H'$ of order $p$ then, as extraspecial groups have order $p^{k}$ for some odd $k$, we have $|\Phi(H)| = p^2$ and a chain $1\trianglelefteq Z(H)\trianglelefteq \Phi(H)\trianglelefteq C_H(\Phi(H))$ satisfying the conditions of Lemma \ref{Oliver1.7}, which again gives a contradiction. Hence, $|Z(H)| = p^2$, and so
$H\leq S$ has exponent $p$, nilpotency class 2, and $|Z(H)| = p^2$. Therefore, $H = \langle x_1, x_2, x_3, x_4\rangle$ such that $x_i^p=1$, $Z(H) = \langle x_3, x_4\rangle$, and $H' = \langle [x_1, x_2]\rangle \leq Z(H)$ has order $p$. Thus, $\langle x_1, x_2\rangle\cong p^{1+2}_+$, which commutes with the remaining generator, hence $H \cong p^{1+2}_+\times C_p$ as claimed.
\end{proof}

In general, there may also be nonabelian $\FF$-essential candidates of order $p^4$ which are isomorphic to $p^{1+2}_+\circ C_{p^2}$. 
We can now determine the candidates for $\FF$-essential subgroups in $S$.

\begin{lemma}\label{EssentialsSUSL} 
If $E\leq S$ is $\FF$-essential then $E\in \{Q, V\}$.
\end{lemma}

\begin{proof}Suppose $E$ is $\FF$-essential but $E\notin \{Q,V\}$.
If  $|E|\leq p^3$ then $|N_S(E)/E| = p$ by \cite[Theorem 1.10]{ValentinaPaperPearls}, which implies  $|E^S| = |S:N_S(E)|\geq p^2$. Notice that as $Z\leq E$ and $Z_2(S) = S'$, we have 
$$EZ_2(S)\leq N_S(E)\trianglelefteq S,$$ 
so every $F\in E^S$ has $F\leq N_S(E)$, and $|E : E\cap Z_2(S)| = p$.
If $|E| = p^2$ then let $K = Z$, and if $|E| = p^3$ then let $K = E\cap Z_2(S)$, of order $p^2$. As $[E\cap Z_2(S), S]\leq Z\leq E\cap Z_2(S)$, we have $K\trianglelefteq S$ in both cases. Then, $K < E < N_S(E)$ with $|N_S(E):K| = p^2$ and $K\trianglelefteq S$, so $|E^S|\leq p+1$, a contradiction. 

If $|E| = p^4$, as $E\neq V$, then  $E$ is nonabelian by Lemma \ref{LUPropertiesOfS}~(5), 
so by Proposition \ref{EssentialCandidatesOrderp4} we have $E\cong p^{1+2}_+\times C_p$, which has rank $3$. Thus, $\Aut(E)$ embeds into $GL_3(p)$ and \cite[Theorem 1.10]{ValentinaPaperPearls} implies that $N_S(E)$ is maximal in $S$. If $E'\neq Z$ then, as $E'< Z(E)$, we have $Z(E) = ZE'$, so $Z(N_S(E)) = Z(E)$. Therefore, $N_S(E) = C_S(Z(E)) > V$, so $N_S(E)\in \XX$ by Lemma \ref{LUPropertiesOfS}~(7). Lemma \ref{LUPropertiesOfS}~(6) then implies that $N_S(E)'\neq Z(N_S(E))$, hence $E'\leq N_S(E)'\cap Z(N_S(E)) = Z$, a contradiction.
Thus, we have $E'=Z$, so that $Z_2(S)$ centralizes the chain $1\leq Z\leq E$, which, as $O_p(\Out_\FF(E)) = 1$, implies that $S'=Z_2(S)\leq E$; therefore, $E\trianglelefteq S$, a contradiction.

 The remaining subgroups are maximal in $S$. Let $M = E$ be a maximal subgroup of $S$. 
  If $M\in \XX$ then  $Z(M)$ and $M'$ are distinct subgroups of $M$ of order $p^2$. We therefore have a chain $M'\trianglelefteq Z(M)M'\trianglelefteq C_M(M')\trianglelefteq M$ of characteristic subgroups of $M$ with successive indices $p$, and Lemma \ref{Oliver1.7} yields a contradiction again.

It only  remains to consider candidates $M\notin \XX$, that is, with $V\nleq M$. Then, Lemma \ref{LUPropertiesOfS}~(7) implies that  $Z(M) = Z(S)$.
If $\Phi(M) = Z$ then $M'=Z$, hence $M=Q$ by Lemma \ref{LUPropertiesOfS}~(4). 
Assume  $M\neq Q$. Since $Q\cap M$ has order $p^4$, we have $(Q\cap M)' = Z$, so $Z\leq M'\leq  \Phi(M)$. 
Thus, any remaining maximal subgroup has $\Phi(M)> Z$. If $\Phi(M) = S'$ then $S$ acts trivially on $M/S'$, contradicting 
Lemma \ref{Oliver1.7}. 

Hence we may assume that $\Phi(M)\neq S'$. Now, by Lemma \ref{LUPropertiesOfS}~(3), $S$  has exponent $p$, so $\Phi(S) = S'$ and $\Phi(M)= M'$. 
%By Lemma \ref{LUPropertiesOfS}~(2), $S$ has nilpotency class $3$, so as $p\geq 5$, $S$ is  regular by \cite[III.10.2 Satz]{Huppert}. Thus, $S^p$
%Now as $Q$ is extraspecial
Thus, $Z<\Phi(M)<S'$, and $|\Phi(M)| = p^2$.  
Note that $Z_2(M)$ has index at least $p^2$ in $M$, so $ S' =Z_2(S) = Z_2(M)$. 
Furthermore, since $\Phi(M)\leq S' = Z_2(M)$, we have $Z_2(M)\leq C_M(\Phi(M))$.
We can therefore build a chain $\Phi(M)\trianglelefteq Z_2(M)\trianglelefteq C_M(\Phi(M)) \trianglelefteq M$ each with index $p$ in the next one, contradicting Lemma \ref{Oliver1.7}. We have now ruled out all subgroups other than $Q$ and $V$, as claimed.
\end{proof}

Now that we know the candidates for $\FF$-essential subgroups, we determine their respective  $\FF$-automizers and the action of these $\FF$-automizers on the relevant modules. We consider $Q$ first. 

\begin{lemma}\label{OutQSL2}Suppose $Q$ is $\FF$-essential. Then $O^{p'}(\Out_\FF(Q))\cong SL_2(p)$ and $Q/Z$ is a direct sum of two natural $SL_2(p)$-modules.

\end{lemma}

\begin{proof} If $Q$ is $\FF$-essential then $O_p(\Out_\FF(Q)) = 1$, whence Lemma \ref{Oliver1.7} implies that $Q/Z$ is a 4-dimensional faithful $\Out_\FF(Q)$-module, where $\Out_\FF(Q)\leq \Out(Q)$ embeds into $GL_4(p)$ with $\Out_S(Q)\cong S/Q$ of order $p$. Since $Z_2(S) = S'$ of order $p^3$ and $[S,S,S] = Z$ by Lemma \ref{LUPropertiesOfS}~(2), we have $ [Q/Z, S]\leq C_{Q/Z}(S) $. Further, since their corresponding subgroups in $S$ are $S'$ and $Z_2(S)$, which coincide, we have $[Q/Z, S] = C_{Q/Z}(S)$. Hence, we can apply Proposition \ref{pNot3V4dQuadratic} to obtain $O^{p'}(\Out_\FF(Q))\cong SL_2(p)$ and $Q/Z$ is a direct sum of two natural $SL_2(p)$-modules.\end{proof}

Now we turn our attention to $V$, which is abelian, so that $\Aut_\FF(V)\cong \Out_\FF(V)$.

\begin{lemma}\label{AutV}Suppose $V$ is $\FF$-essential. Then $O^{p'}(\Aut_\FF(V))\cong PSL_2(p^2)$ and $V$ is a natural $\Omega_4^-(p)$-module for $O^{p'}(\Aut_\FF(V))$.
\end{lemma}

\begin{proof}Suppose $V$ is $\FF$-essential. Then the Sylow $p$-subgroups of $\Aut_\FF(V)$  are elementary abelian of order $p^2$, $\Aut_\FF(V)$ contains a strongly $p$-embedded subgroup and, since $V\cong C_p^4$, it embeds into $\Aut(V)\cong GL_4(p)$. Furthermore, since $\Aut_S(V)$, a Sylow $p$-subgroup of $O^{p'}(\Aut_\FF(V))$, fixes the $1$-subspace $Z\leq V$, Proposition \ref{speinGL4} implies that $O^{p'}(G)\cong PSL_2(p^2)$ and $V$ is a natural $\Omega_4^-(p)$-module.
\end{proof}

We can now determine the relationship between the $\FF$-essential subgroups and $O_p(\FF)$.

\begin{lemma}\label{EssentialsInSU4} Both  $Q$ and $V$ are $\FF$-essential if and only if $O_p(\FF)=1$. 
\end{lemma}
\begin{proof}By Lemma \ref{EssentialsSUSL}, the $\FF$-essential candidates are $Q$ and $V$, both of which are characteristic in $S$, hence normalized by $\Aut_\FF(S)$. If $L\in \{Q, V\}$ is the only $\FF$-essential subgroup then the Alperin-Goldschmidt Fusion Theorem \ref{Alperin} shows that  $\FF = \langle \Aut_\FF(S), \Aut_\FF(L)\rangle$, where $L$ is normalized by $\Aut_\FF(S)$ and $\Aut_\FF(L)$, hence %by Theorem \ref{NormalSubgroup} 
we have $L \trianglelefteq \FF$ by Proposition \ref{AKONormal}, which contradicts $O_p(\FF) = 1$. Thus, if $O_p(\FF) = 1$, then both $V$ and $Q$ are $\FF$-essential.

Conversely, assume both $Q$ and $V$ are $\FF$-essential. By Proposition \ref{AKONormal}, $O_p(\FF)\leq Q\cap V$ is a subgroup of $V$. But by Lemma \ref{AutV} $V$ is a natural $\Omega_4^-(p)$-module for $O^{p'}(\Aut_\FF(V))$, in particular $O^{p'}(\Aut_\FF(V))$ acts irreducibly on $V$, hence $O_p(\FF) = 1$.
\end{proof}

We can now prove Theorem \ref{essentials}.
\begin{proof}[Proof of Theorem \ref{essentials}]
Follows from Lemmas \ref{EssentialsSUSL}, \ref{OutQSL2}, \ref{AutV} and  \ref{EssentialsInSU4}, .
\end{proof}

\section{Classification of the fusion systems on \texorpdfstring{$S$}{\textit{S}}}\label{sec:classification}

We continue with the assumptions that $S$ is a Sylow $p$-subgroup of $SU_4(p)$ and $\FF$ is a saturated fusion system on $S$. 
We now have all the pieces required to construct $\FF$, it remains to put them together. 
We begin with the case when $O_p(\FF) \neq 1$, where we consider $p\geq 3$.

\begin{prop}\label{OpNotTrivial}Suppose $p$ is odd and $O_p(\FF)\neq 1$. Then $\FF$ has a model. In particular, $\FF$ is realizable by a finite group.
\end{prop}
\begin{proof} If $\FF$ has no $\FF$-essential subgroups, it is realizable by $S\rtimes \Out_\FF(S)$. 
If $p\geq 5$ then Lemma \ref{EssentialsInSU4} implies that $\FF$ has at most one $\FF$-essential subgroup. 
If $p=3$, the only $\FF$-essential candidates are $Q$ and $V$ by \cite[Proposition 7]{BacMcL}, and if both $Q$ and $V$ are essential then  \cite[Proposition 12~(iii)]{BacMcL} implies that $O_3(\FF) = 1$. Thus, whenever  $O_3(\FF) \neq 1$, $\FF$ contains at most one $\FF$-essential subgroup. 
In either case, let $L\in \{Q, V\}$ be the only $\FF$-essential subgroup. As $L$ is characteristic in $S$, Proposition \ref{AKONormal} implies that $L\trianglelefteq \FF$, and $L$ is $\FF$-centric as it is $\FF$-essential. Therefore, $\FF$ is constrained, and the Model Theorem for constrained fusion systems \cite[Theorem I.4.5]{AKO} implies that $\FF$ has a model, that is a finite group realizing $\FF$.
\end{proof}

We assume that $p\geq 5$ and $O_p(\FF) = 1$ for the rest of this section. In particular, the $\FF$-essential subgroups are $Q$ and $V$ by Lemma \ref{EssentialsInSU4}.
The goal of this section is to complete the proof of Theorem \ref{CaseU} by proving the following result.

\begin{thm}\label{ClassificationCaseU} Assume $p\geq 5$ and that $S$ is a Sylow $p$-subgroup of $SU_4(p)$. Then there is a one-to-one correspondence between saturated fusion systems $\FF$ on $S$ with $O_p(\FF)=1$ and groups $G$ with $PSU_4(p)\leq G \leq \Aut(PSU_4(p))$ which realize them. In particular, there are no exotic fusion systems $\FF$ on $S$ with $O_p(\FF)=1$. There are $5$ such $\FF$ up to isomorphism when $p\equiv 1\pmod 4$ and $8$ when $p\equiv 3\pmod 4$.
\end{thm}

The main part of the proof is dedicated to determining uniqueness of the fusion subsystem $\FF_0$, which we prove in Proposition \ref{UniquenessOutQInSU}. 
We begin by setting up some notation for maps that we will use throughout. Recall Definition \ref{AutSFromE} where, for each $\FF$-essential subgroup $E$, we set
$$\Aut_\FF^E(S) := \langle \alpha \in \Aut_\FF(S) \mid \alpha|_E\in O^{p'}(\Aut_\FF(E))\rangle.$$

\begin{notation}\label{Notation5}
\begin{itemize} \hspace{2em}
\item $D_V$ is a complement to $\Aut_S(V)$ in  $N_{O^{p'}(\Aut_\FF(V))}(\Aut_S(V))$.

\item $D_Q$ is a complement to $\Aut_S(Q)$ in  $N_{O^{p'}(\Aut_\FF(Q))}(\Aut_S(Q))$.

\item $L_V$ is a complement to $\Inn(S)$ in $\Aut_\FF^V(S)$, with $L_V|_V = D_V$.

\item $L_Q$ is a complement to $\Inn(S)$ in $\Aut_\FF^Q(S)$, with  $L_Q|_Q = D_Q$.

\item $d = \gcd(4,p+1) = |Z(SU_4(p))|$.  %$= |SU_4(p)|/|PSU_4(p)|$.
\end{itemize}
\end{notation}

We remark that, by Lemma \ref{OutFSBijection}, $D_V$ and $L_V$ exist and satisfy $D_V\cong L_V$ and $D_Q\cong L_Q$, and
that they depend on the choice of $\FF$. Furthermore, since $V$ and $Q$ are characteristic in $S$, $\Aut_\FF^V(S)$ and $\Aut_\FF^Q(S)$ normalize each other, so that $\Aut_\FF^0(S) = \Inn(S)L_VL_Q$ by definition.

%Our strategy in this section is the following:
%\begin{enumerate}
%\item For each arbitrary conjugate of $O^{p'}(\Aut_\FF(V))\cong PSL_2(p^2)$ in $\Aut(V)$ which contains $\Aut_S(V)$, we study $D_V$, $L_V$, and $C_{L_V}(Z)$.
%\item For each arbitrary conjugate of $O^{p'}(\Aut_\FF(Q))\cong SL_2(p)$ in $\Aut(Q)$ which contains $\Aut_S(Q)$, we study $D_Q$ and $L_Q$.
%\item In particular, we determine the conjugacy class in $Sp_4(p)\leq \Out(Q)$ of the actions induced by generators of  $C_{L_V}(Z)$ and $L_Q$.
%
%\item We study $L_V\cap L_Q$, which allows us to determine $\Aut_\FF^0(S)$ uniquely. 
%
%
%\item We fix a copy of $\Aut_\FF^0(S) = \Aut_{\FF_0}(S)$ and use the subgroup $T$ of morphisms centralizing both $Z$ and $S/Q$  
%to determine $O^{p'}(\Aut_{\FF_{0}}(V))$ and $O^{p'}(\Aut_{\FF_{0}}(Q))$ uniquely. This determines $\FF_0$ uniquely up to isomorphism of fusion systems. 
%
%\item We realize $\FF_0$ as $\FF_S(PSU_4(p))$, proving that $\FF_0$ is saturated.
%\item Finally, we obtain a one-to-one correspondence between the saturated fusion systems on $S$ with $O_p(\FF) = 1$ and the intermediate subgroups between $PSU_4(p)$ and $\Aut(PSU_4(p))$
%and count them, concluding the proof of Theorems \ref{ClassificationCaseU}, \ref{CaseU} and \ref{mainthm}.
%
%\end{enumerate}
%

We now begin the proof. We first obtain the action on $V$ of $D_V$ and $L_V$. The following lemma uses Lemma \ref{Omega4-Calculation} to set notation.

\begin{lemma}\label{LiftsFromV} $D_V$ is cyclic of order $(p^2-1)/2$, and we have: 
\begin{enumerate}
\item $D_V$ acts on $Z$ and $V/S'$ as an element of order $p-1$, and $S'/Z$ as an element of order $p+1$. Furthermore, all the actions are irreducible. 

\item $L_V$ has order $(p^2-1)/2$ %$\widetilde t_V|_V = t_V$, 
and a generator of $L_V|_Q\leq  \Aut_\FF(Q)$ acts on $Z$ and $S/Q$ as an automorphism of order $p-1$.

\item Let $\sigma\in L_V$ be the unique involution. When $4\mid p+1$, $\sigma$ centralizes  $Z$ and $V/S'$ and inverts $S'/Z$  whereas, when $4\mid p-1$, $\sigma$ centralizes $S'/Z$ and inverts $Z$ and $V/S'$. In both cases it inverts $S/V$.
\end{enumerate}
\end{lemma}

\begin{proof}By Theorem \ref{essentials}~(3), $V$ is a natural $\Omega_4^-(p)$-module for $O^{p'}(\Aut_\FF(V))\cong PSL_2(p^2)$,  
and the structure of this module was described 
in Lemma   \ref{Omega4-Calculation}, where a generator of $D_V$ corresponds to $t$. We extract the notation from there. In particular, Lemma \ref{Omega4-Calculation}~(2) implies that 
$D_V$ acts on $V$ in the way described in part (1). 
Recall that $S = VQ$, $V\cap Q = Z_2(S) = S'$ and thus $S/Q\cong V/S'$.
Then $Z = C_V(S) = C_V(\Aut_S(V)) = V_3$ and $[V, S] = [V, \Aut_S(V)] = S' = V_2\oplus V_3$ as in Lemma \ref{Omega4-Calculation}~(2), 
whence  part (2) follows by Lemma \ref{OutFSBijection}. 
Finally, for part (3), the action of $\sigma$ on $V$ and $S/V\cong \Aut_S(V)$ is described in Lemma \ref{Omega4-Calculation}~(4). 
\end{proof}

Using similar arguments, we obtain $C_{L_V}(Z)$, which will help us determine $O^{p'}(\Out_\FF(Q))$ uniquely as a subgroup of $Sp_4(p)$ (not just up to conjugacy). In the next two lemmas we exploit \cite{SrinivasanSp4}, for which we require the eigenvalues of the action of generators of $L_V$ and $L_Q$ on $Q$. The condition that $p\geq 5$ is used in the next lemma.

\begin{lemma}\label{ConjClassLift} Assume $p\geq 5$, then
 $C_{L_V}(Z) \cong C_{(p+1)/2}$ centralizes both $Z$ and $S/Q$. 
Let $\psi$ be a generator of $C_{L_V}(Z)$, then $\psi |_Q\Inn(Q)
\in C_{\Out(Q)}(Z)\cong Sp_4(p)$ is in the $Sp_4(p)$-conjugacy class $B_6(2)$ of \cite{SrinivasanSp4}, with $|C_{Sp_4(p)}(\psi|_Q\Inn(Q))| = |GU_2(p)|$.
\end{lemma}

\begin{proof}

By Lemma \ref{Omega4-Calculation}~(3), we have $C_{D_V}(C_V(\Aut_S(V))) = C_{D_V}(Z)\cong C_{(p+1)/2}$. Let $\psi$ be a generator of $C_{L_V}(Z)$. By Lemma \ref{LiftsFromV}~(2), $C_{L_V}(Z)\cong C_{(p+1)/2}$, and $\psi$ centralizes $Z$ and $V/S'$ . 
As $\psi|_V\in D_V\leq O^{p'}(\Aut_{\FF}(V)$, we have $\det(\psi|_V) = 1$ and, as $\psi$ has order $(p+1)/2$,  $p\geq 5$, and $\psi$ centralizes $Z$ and $V/S'$, 
the eigenvalues of $\psi$ on $S'/Z$ are in $GF(p^2)\setminus GF(p)$ and are $\eta^{2}$, $\eta^{-2}$ for $\kappa$ a primitive element of $GF(p^2)$.

Note that the action of $\psi$ needs to be consistent with the commutator structure of $S$. Let $s\in V\setminus S'$, then $(sS')\psi = sS'$, and let $q,r\in Q$ and $x\in S'$.
We consider a homomorphism $\theta : Q\rightarrow S'/Z$ defined by $q\theta = [q, sS']Z$. Then  $$qr\theta = [qr, sS']Z = [q,sS']^r[r,sS']Z = [q,sS'][r,sS']Z = q\theta r\theta,$$ 
hence $\theta$ is a homomorphism.
Further, $\theta$ preserves the action of $\psi$ as
$$q\theta \psi = [q,sS']Z\psi = [q\psi,(sS')\psi]Z = [q\psi,sS']Z = q\psi\theta$$ since $\psi$ centralizes $sS'$ and $Z$. Since $\ker\theta = C_Q(V) = S'$, we conclude that $Q/S'$ and $S'/Z$ are isomorphic as $\langle \psi\rangle$-modules and we obtain the eigenvalues of the projection of $\psi$ to $Q/Z$, which are $\eta^2$, $\eta^2$, $\eta^{-2}$ and $\eta^{-2}$.
Recall that $C_{\Out(Q)}(Z)\cong Sp_4(p)$ by Theorem \ref{AutomorphismOfExtraspecial}.
 This determines the conjugacy class of $\psi|_Q \Inn(Q)$ 
as an element of $C_{\Out(Q)}(Z)\cong Sp_4(p)$ to be the class $B_6(2)$ in the notation of \cite{SrinivasanSp4},
which has $|C_{Sp_{4}(p)}(\psi)| = p(p+1)(p^2-1) = |GU_2(p)|$. 
\end{proof}

We now consider $D_Q$ and $L_Q$.

\begin{lemma}\label{LiftsFromQU} $L_Q$ is cyclic of order $p-1$ and centralizes $Z$.
Let $\tau$ be a generator of $L_Q$, then $\tau|_Q\Inn(Q)$ is in the $Sp_4(p)$-conjugacy class $B_3(1,1)$ of \cite{SrinivasanSp4}.
\end{lemma}
\begin{proof}By Theorem \ref{essentials}~(2), $O^{p'}(\Out_\FF(Q))\cong SL_2(p)$ and acts on $Q/Z$ as a direct sum of two natural $SL_2(p)$-modules. 
By Lemma \ref{OutFSBijection}, we have
$$L_Q\cong \Aut_\FF^Q(S)/\Inn(S)\cong N_{O^{p'}(\Aut_\FF(Q))}(\Aut_S(Q))/\Aut_S(Q)\cong D_Q\cong  C_{p-1}.$$
Now $\tau|_Q\in N_{O^{p'}(\Aut_\FF(Q))}(\Aut_S(Q))$ by definition, so $\tau$ centralizes $Z$ and so does $L_Q$. Finally, 
as in the proof of Proposition \ref{pNot3V4dQuadratic}, we see that the projection of $\tau|_Q$ 
to $Q/Z$ has eigenvalues $\lambda, \lambda, \lambda^{-1}$ and $\lambda^{-1}$, where $\lambda$ is a primitive element of $GF(p)$. Therefore, $\tau|_Q\Inn(Q)\in C_{\Out(Q)}(Z)\cong Sp_4(p)$ 
is in the conjugacy class denoted by $B_3(1,1)$ in \cite{SrinivasanSp4}.
\end{proof}

We now calculate the intersection of $L_V$ and $L_Q$ to obtain $\Aut_\FF^0(S)$ which, together with $O^{p'}(\Aut_\FF(Q))$ and $O^{p'}(\Aut_\FF(V))$, will generate the fusion system $\FF_0$ on $S$.
 
\begin{prop}\label{AllLiftsU} The group $\Aut_\FF^0(S) = \Inn(S)L_QL_V$ is uniquely determined as a subgroup of $\Aut_{\Aut(SU_4(p))}(S)$.
In particular, $\Aut_\FF^0(S)\cong \Inn(S)\rtimes(C_{p-1}\circ_\Delta C_{(p^2-1)/2}),$ where $\Delta = L_Q\cap L_V$ has order $d/2$. 
\end{prop} 
\begin{proof}

Since the only $\FF$-essential subgroups are $Q$ and $V$, $\Aut_\FF^0(S)$ is generated by 
$\Inn(S)$, $L_Q$ and $L_V$. By Lemma \ref{LiftsFromQU}, we have $L_Q = C_{L_Q}(Z)\cong C_{p-1}$, and
Lemma \ref{LiftsFromV}~(2) yields $L_V\cong C_{(p^2-1)/2}$ with $|L_V:C_ {L_V}(Z)| = p-1$.
It remains to consider $L_Q\cap L_V\leq C_{L_V}(Z)$, which has order at most $(p+1)/2$. As $L_Q$ has order ${p-1}$, the intersection can have order at most $\gcd(p-1, (p+1)/2) \leq 2$. Therefore, if $4\nmid p+1$, that is $d=2$, the intersection is trivial. 

However, if $4\mid p+1$, that is $d=4$, we can have intersection of size at most $2$. 
In Lemma \ref{LiftsFromQU} we proved that the action induced by $\tau$ on $Q/Z$ is in $Sp_4(p)$-conjugacy class $B_3(1,1)$, thus 
$\tau^{(p-1)/2}$ acts on $Q/Z$ as $-I_4$, and it centralizes $Z$ and $S/Q$.
Similarly, as $4\mid p+1$, Lemma \ref{LiftsFromV}~(3) implies that $\sigma$ centralizes $Z$ and $S/Q$ and inverts $Q/Z$, in other words, $\sigma = \tau^{(p-1)/2}$. 

Hence, the order of $\Aut_\FF^0(S)$ is as claimed. Finally, Lemma \ref{UniquenessOfOutF0S} implies that there is a unique subgroup of $\Aut_{\Aut(SU_4(p))}(S)$ of the shape described in either case. In particular, the isomorphism type of $\Aut_\FF^0(S)$ is determined.
\end{proof}

We can now use the fact that $\Aut_\FF^0(S)$ is uniquely determined to show that this determines $O^{p'}(\Aut_\FF(V))$ and $O^{p'}(\Aut_\FF(Q))$ uniquely, which proves uniqueness of the subsystem $\FF_0$.

\begin{prop}\label{UniquenessOutQInSU} $\FF_0$ is unique up to isomorphism. 
\end{prop}
\begin{proof}
We have determined $\Aut_{\FF}^0(S)\cong \Inn(S)\rtimes (C_{p-1}\circ_\Delta C_{(p^2-1)/2})$ uniquely as a subgroup of $\Aut_{PSU_4(p)}(S)$ in Proposition \ref{AllLiftsU}. Fix this subgroup. 
Then, by Lemma \ref{Op'AndSGenerate},
we need to determine $O^{p'}(\Aut_{\FF_0}(Q))$ and $O^{p'}(\Aut_{\FF_0}(V))$ uniquely.

We now consider the subgroup 
$$T = \{\alpha\in \Aut_{\FF}^0(S)\mid z\alpha = z\text { for all } z \in Z \text{ and } [S,\alpha]\le Q\} \leq \Aut_{\FF_0}(S).$$
Then, $C_{L_V}(Z)\leq T$ by Lemma \ref{ConjClassLift}, and $T\cap L_Q$ has order $2$ as the second condition is equivalent to centralizing $S/Q\cong \Aut_S(Q)$.
Thus, $T$ has index $(p-1)^2/2$ in $\Aut_{\FF_0}(S)$ and shape $\Inn(S) \rtimes (C_2\circ_\Delta C_{(p+1)/2})$
In particular, a complement to $\Aut_S(V)$ in $T$ is cyclic of order ${2(p+1)/d}$.

We now consider the subgroup $T_V$ of $\Aut_{\FF_0}(V)$ obtained by restricting maps in $T$ to $V$. This yields $T_V\leq N_{\Aut_{\FF_0}(V)}(\Aut_S(V)) $ of shape  $\Aut_S(V)\rtimes C_{2(p+1)/d}$, which in turn determines $N_{\Aut_{\FF_0}(V)}(\Aut_S(V))$ uniquely. Then,  Lemma \ref{Omega4-Calculation}~(5) implies that there is a unique non-degenerate symplectic form preserved by $N_{\Aut_{\FF_0}(V)}(\Aut_S(V))$, which determines $O^{p'}(\Aut_{\FF_0}(V))$ uniquely.

We now turn our attention to $Q$.  We have $O^{p'}(\Out_{\FF_0}(Q))\cong SL_2(p)$ by Theorem \ref{essentials}~(2). 
There is some  $\beta\in T$ of order $(p+1)/2$ satisfying $\beta|_V = t_V^{p-1}$, then $\beta$ acts on $S$ like $\psi$. Then, by Lemma \ref{ConjClassLift}, 
$\beta|_Q\Inn(Q)$ acts on $Q/Z$ via a matrix in $Sp_4(p)$-conjugacy class $B_6(2)$, 
which satisfies $|C_{Sp_4(p)}(\beta|_Q\Inn(Q))| = |GU_2(p)|$.

Since $O^{p'}(\Aut_{\FF_0}(Q))\trianglelefteq \Aut_{\FF_0}(Q)$, $\beta$ normalizes $O^{p'}(\Aut_{\FF_0}(Q))$. 
We have $\Aut(SL_2(p))\cong PGL_2(p)$ and if $R\in \Syl_p(PGL_2(p))$ then $R = C_{PGL_2(p)}(R)$. Thus, as $\Aut_S(Q)/\Inn(Q)$ is centralised by $\beta$, $\beta$ centralizes $O^{p'}(\Out_{\FF_0}(Q))$.
That is, $O^{p'}(\Out_\FF(Q))\leq C_{Sp_4(p)}(\beta|_Q\Inn(Q))$, which contains a unique subgroup isomorphic to $SL_2(p)$. Since $\Aut_{\FF}^0(S)$ also contains a subgroup acting transitively on $Z$, we conclude that $O^{p'}(\Aut_{\FF_0}(Q))$ is uniquely determined as a subgroup of $\Aut(Q)$.

As $\FF_0 = \langle O^{p'}(\Aut_\FF(V)), O^{p'}(\Aut_\FF(Q)), \Aut_\FF^0(S)\rangle$ by definition,
we have shown that fixing $\Aut_\FF^0(S)$ uniquely determines $\FF_0$. 
\end{proof}

At this stage we determine the fusion system of $PSU_4(p)$. Recall that the fusion systems of $SU_4(p)$ and $PSU_4(p)$ are isomorphic by Lemma \ref{IsomorphismQuotient}. 

\begin{lemma}\label{FusionU4} $\FF_S(PSU_4(p))$ is isomorphic to $\FF_0$ whenever $p\geq 5$. 
In particular, $\FF_0$ is saturated and any saturated fusion system $\FF$ on $S$ with $O_p(\FF) = 1$ satisfies $O^{p'}(\FF)\cong \FF_0$.
\end{lemma}

\begin{proof} Let $\GG := \FF_S(PSU_4(p))$, which is a saturated fusion system.
By \cite[Table 8.10]{BHRD}, there are two conjugacy classes of maximal parabolic subgroups of $SU_4(p)$, with structure
$N_{SU_4(p)}(Q) \sim p^{1+4}:SU_2(p):C_{p^2 -1}$ and  $N_{SU_4(p)}(V)\sim C_p^4:SL_2(p^2):C_{p - 1}$, whose intersection is the Borel subgroup $N_{SU_4(p)}(V)$ of shape $S : (C_{p-1}\times C_{p^2-1})$.  Note that these all contain $Z(SU_4(p))$, which has order $d$ by \cite[Table 2.1.D]{KLClassical}. In particular, $N_{SU_4(p)}(V)$ acts irreducibly on $V$, which proves that $O_p(\GG) = 1$. Finally, for $E\in \{S,Q,V\}$, we have $\Aut_{\GG}(E)\cong \Aut_{\FF_0}(E)$, whence Proposition \ref{UniquenessOutQInSU} implies that $\GG$ is isomorphic to $\FF_0$.  Therefore, $\FF_0$ is saturated and, since $\FF_0$ is uniquely determined, Theorem \ref{PrimeIndexFusion} implies that any saturated $\FF$ on $S$ with $O_p(\FF) =1$ satisfies $O^{p'}(\FF) \cong \FF_0$.
\end{proof}

We have now constructed and realized the unique smallest possible fusion system $\FF_0$ with $O_p(\FF)=1$ on $S$. Recall that any saturated fusion system on $S$ with $O_p(\FF) = 1$ satisfies $\FF = \langle \FF_0, \Aut_\FF(S)\rangle$ by Lemma \ref{Op'AndSGenerate}. We now consider the largest possible $\Aut_\FF(S)$, which by Lemma \ref{OutS} has $|\Aut_\FF(S)| = p^52(p+1)(p-1)^2$, to conclude the proof of Theorem \ref{ClassificationCaseU}.

\begin{proof}[Proof of Theorem \ref{ClassificationCaseU}]
Let $\HH := \FF_{S}(\Aut(PSU_4(p)))$, then $\HH$ is a saturated fusion system on $S$ containing $\GG = \FF_S(PSU_4(p))$, hence $O_p(\HH) = 1$. Furthermore, Lemma \ref{OutS} implies that $\Aut_{\HH}(S)$ is as large as possible. Therefore, any saturated fusion system $\FF$ on $S$ with $O_p(\FF) = 1$ is isomorphic to a fusion system $\EE$ intermediate between $\FF_0\cong \GG$ and $\HH$, which are in a one-to-one correspondence with subgroups satisfying $\Aut_{\FF_0}(S)\leq \Aut_\EE(S)\leq \Aut_{\HH}(S)$
by Theorem \ref{PrimeIndexFusion},
and each of them is realized by a corresponding intermediate subgroup between $PSU_4(p)$ and $\Aut(PSU_4(p))$.

Finally, the number of non-isomorphic fusion systems corresponds to the number of conjugacy classes of subgroups of $\Aut(PSU_4(p))/PSU_4(p)$ under $\Aut(S)$. 
When $p\equiv 1\pmod 4$, we have $\Aut(PSU_4(p))/PSU_4(p)\cong C_2\times C_2$, which has $5$ conjugacy classes of subgroups, whereas when $p\equiv 3\pmod 4$, we have $\Aut(PSU_4(p))/PSU_4(p)\cong D_8$, which has $8$ conjugacy classes of subgroups, concluding our claims. 
Any further isomorphism between these fusion systems would be in $\Aut(S)$, hence it would be a $p$-element by Lemma \ref{OutS} which is impossible as $p\geq 5$, whence there can be no more isomorphisms between the fusion systems enumerated above.
\end{proof}

\section*{Acknowledgements}

The main results of this paper were proved by the author in his PhD Thesis, 
%at the School of Mathematics, University of Birmingham, United Kingdom, 
under the supervision of Professor Chris Parker. The author is extremely grateful to him for his support. This work was supported by the School of Mathematics at the University of Birmingham.

%This research did not receive any specific grant from funding agencies in the public, commercial, or not-for-profit sectors.

Declarations of interests: none.

\newcommand{\etalchar}[1]{$^{#1}$}
\providecommand{\bysame}{\leavevmode\hbox to3em{\hrulefill}\thinspace}
\providecommand{\MR}{\relax\ifhmode\unskip\space\fi MR }
\providecommand{\MRhref}[2]{%
  \href{http://www.ams.org/mathscinet-getitem?mr=#1}{#2}
}
\providecommand{\href}[2]{#2}

\end{document}